\RequirePackage{fix-cm}
\documentclass[smallcondensed]{svjour3}     
\usepackage{graphicx}

 \usepackage{mathptmx}      
\usepackage{amsmath}
\usepackage{amsfonts}
\usepackage{caption,subcaption}
\usepackage{setspace}
\usepackage{enumitem}
\usepackage{textcomp}

\makeatletter
\renewcommand\section{\@startsection{section}{1}{\z@}%
                                   {-3.5ex \@plus -1ex \@minus -.2ex}%
                                   {2.3ex \@plus.2ex}%
                                   {\normalfont\large\bfseries}}
\makeatother

\usepackage[top = 1in, left = 1in, right = 1in, bottom = 1in]{geometry}

\newenvironment{myfont}{\fontfamily{phv}\selectfont}{\par}
\def\b{\begin{myfont}}
\def\t{\end{myfont}}
\def\D{\mbox{\b D\t}}
\def\d{\mbox{\b d\t}}
\def\E{\mbox{\b E\t}}
\def\e{\mbox{\b e\t}}

      \newcommand{\Z}{{\mathbb Z}}

      \newcommand*{\defeq}{\mathrel{\vcenter{\baselineskip0.5ex \lineskiplimit0pt
                     \hbox{\scriptsize.}\hbox{\scriptsize.}}}%
                     =}
\newcommand*{\eqdef}{\mathrel=\hspace{-.041in}{\vcenter{\baselineskip0.5ex \lineskiplimit0pt
                     \hbox{\scriptsize.}\hbox{\scriptsize.}}}%
                     \hspace{.04in}}

\newcommand{\UT}{a}             
\newcommand{\firstcompleted}{b}
\newcommand{\secondcompleted}{c}
\newcommand{\thirdcompleted}{d}
\newcommand{\diamondtoclassic}{e}
\newcommand{\secondDAI}{f}
\newcommand{\DAIUC}{g}
\newcommand{\NCCUC}{h}
\newcommand{\MaxsDistinct}{i}
\newcommand{\FT}{j}
\newcommand{\classiccombo}{k}
\newcommand{\MaxsDistinctMaxsFree}{l}
\newcommand{\secondNCCMaxsFree}{m}
\newcommand{\fullclassictodiamond}{n}
\newcommand{\UPUE}{o}

\newcommand{\UMCLMEE}{a}
\newcommand{\RankCLE}{b}
\newcommand{\Ss}{c}
\newcommand{\firstNTC}{d}
\newcommand{\classictodiamond}{e}
\newcommand{\thirdNCCsecondNTCUC}{f}

\newcommand{\firstNODk}{a}
\newcommand{\generalclassiccombo}{b}
\newcommand{\secondNODk}{c}
\newcommand{\firstDkMD}{d}
\newcommand{\secondDkMD}{e}
\newcommand{\UCkDkMF}{f}

\numberwithin{theorem}{section}

\spnewtheorem{fact}[theorem]{Fact}{\bf}{\it}

\spnewtheorem*{axiom*}{}{}{\it}

\spnewtheorem*{property*}{}{}{\it}

\spnewtheorem{conj}[theorem]{Conjecture}{\bf}{\it}

\spnewtheorem{myremark}[theorem]{Remark}{\bf}{\rm}
\spnewtheorem*{myremark*}{Remark}{\bf}{\rm}
\spnewtheorem*{Added Notes}{Added Notes}{\bf}{\rm}

\spnewtheorem{mylemma}[theorem]{Lemma}{\bfseries}{\itshape}

\spnewtheorem{myproposition}[theorem]{Proposition}{\bfseries}{\itshape}

\spnewtheorem{mycorollary}[theorem]{Corollary}{\bfseries}{\itshape}
\spnewtheorem{mydefinition}[theorem]{Definition}{\bfseries}{\itshape}

\numberwithin{table}{section}
\numberwithin{figure}{section}

\begin{document}

\large

\title{$\bold{\emph{d}}$-Complete Posets:  \\ Local Structural Axioms, Properties, and Equivalent Definitions
 \thanks{This paper is based upon an April 2015 UNC Chapel Hill Masters project by Lindsey M. Scoppetta that was cowritten with project advisor Robert A. Proctor.}
}

\titlerunning{$d$-Complete Posets: Local Structural Axioms, Properties, and Equivalent Definitions}        

\author{Robert A. Proctor         \and
        Lindsey M. Scoppetta 
}


\institute{Robert A. Proctor \at
              Department of Mathematics, University of North Carolina at Chapel Hill \\
              Tel.: 919-962-9623\\
              \email{rap@email.unc.edu}           
           \and
           Lindsey M. Scoppetta \at
              \email{lindseyscoppetta@gmail.com}           
}

\date{Preprint date: \today }

\maketitle

\begin{abstract}
Although  $d$-complete posets arose along the interface between algebraic 
combinatorics and Lie theory,  they are defined using only requirements on 
their local structure.
These posets are a mutual generalization of rooted trees, shapes, and 
shifted shapes.
They possess Stanley's hook product property for their $P$-partition generating functions and Sch\"{u}tzenberger's well defined jeu de taquin rectification property.
The original definition of  $d$-complete poset was lengthy,  but 
more succinct definitions were later developed.
Here several definitions are shown to be equivalent.
The basic properties of $d$-complete posets are summarized.
Background and a partial bibliography for these posets is given.

\keywords{$d$-complete poset \and $\lambda$-minuscule element \and double tailed diamond \and hook length poset \and partially ordered set}
\end{abstract}


\section{Introduction}

\label{intro}
This paper is entirely poset-theoretic,  apart from some 
background and motivational remarks.
$d$-Complete posets arose in the area of overlap between combinatorics, 
representation theory, and algebraic geometry that is inhabited by Young 
tableaux, Coxeter (Weyl) groups, Kac-Moody Lie algebras, and flag varieties.
Generalizing (shifted) Young diagrams,  $d$-complete posets have been shown to possess 
both Stanley's hook product property for their  $P$-partition generating 
functions \cite{Japan} and Sch\"{u}tzenberger's well defined jeu de taquin rectification 
property \cite{JDT}. Specializing these hook identities gives generalizations of the FRT hook formula for the number of standard Young tableaux to enumerations of the linear extensions of $d$-complete posets.

	$d$-Complete posets can be defined with various combinations of 
local structural axioms.
They have been classified with Dynkin diagrams \cite{DDCT}.
We study the interplay between a number of structural axioms and indicate 
which combinations of these axioms produce some useful local structural 
properties.
Both locally finite (all intervals are finite) and finite posets are 
considered.
Sections~\ref{k=3 propositions} and~\ref{Propositions kgeq3} contain conjectures that graph theorists may be able to confirm by arguing that certain rank size growths are unbounded.
We prove that the lengthy original published \cite{DDCT} 
definition of ``$d$-complete" is equivalent to the succinct most recent \cite{Japan} definition.
This  is the first step needed in a 
sequence of journal papers that will provide a complete 
derivation of the multivariate hook product identity for colored $d$-complete posets \cite{Japan} that 
was obtained with Dale Peterson. 
We are hoping this paper becomes a standard reference for $d$-complete posets:
Given the growing interest in these posets, we summarize their properties, outline the foremost background material, and provide a partial bibliography.

	Before \cite{DDCT} was written,  in 1994 an earlier notion of (colored)
$d$-complete poset was developed while combinatorializing \cite{Mar} a basis 
theorem of Seshadri for certain representations of simple Lie algebras.
There the structural axioms arose from relations within the universal 
enveloping algebra of the Lie algebra.
When the jeu de taquin rectification algorithm was shown to be well defined for 
$d$-complete posets, the remarkable ``simultaneous" property was also 
obtained.
It seemed likely that some algebraic phenomena related to the local 
structure of  $d$-complete posets underlay this property;  these phenomena 
would be related to the reduced decompositions of  $\lambda$-minuscule elements of 
Weyl groups.
It is hoped that the understanding of the interplay among the local 
structural axioms obtained here will facilitate the development of an 
explanation of the interplay between the order theoretic structure of 
$d$-complete posets and some of the algebraic structures in this area of 
mathematics.

	For a full understanding of the algebraic roles of  $d$-complete 
posets,  the notion of colored  $d$-complete poset is needed.
The notion of colored  $d$-complete was shown \cite{Wave} to be essentially 
equivalent to the purely structural notion of  $d$-complete considered here.
This paper sets the stage for a sequel in which the 
lengthy original definition of colored  $d$-complete poset is shown to be 
equivalent to the shorter current \cite{Japan} colored definition.
This equivalence is also needed for the journal papers version of the 
conference proceedings contribution \cite{Japan}. 
To prepare for developing a notion of colored $d$-complete for locally finite posets, 
here we are careful to delineate between the finite and locally finite cases. See Figure~\ref{e7(7) E6}b.

	Ishikawa and Tagawa have developed a category of finite posets 
that vastly extends the category of finite colored $d$-complete posets.
They have shown \cite{IsTa} that their ``leaf'' posets also possess hook product
identities for the associated  colored $P$-partition generating functions; these identities subsume those in \cite{Japan}.
Presently leaf posets are defined only by the presentation of families of 
Hasse diagrams that generalize the families of diagrams appearing in the 
classification \cite{DDCT} of  $d$-complete posets.
The hook property is so special and this extension of it is so nice that 
it is natural to expect that there exists some underlying algebraic or 
geometric explanation for it;  currently their combinatorial 
generating function calculations proceed class-by-class.
As a first step toward a uniform understanding,  it would be desirable for 
someone to develop an axiomatic definition of the notion of leaf poset in 
the spirit of the axiomatic considerations presented here.
For example, the ``short intervals are small" property considered here may be useful for 
studying leaf posets,  as well as for the study of colored  $d$-complete 
posets.

	In the original definition of  $d$-complete poset,  there were three local structural conditions for each positive integer  $k\geq3$.
Although we show that there exist useful alternates to or hybrids of these 
three axioms,  for each  $k\geq3$  each of the following three aspects remains 
present:
First, if a convex subset of the poset is isomorphic to the poset formed by 
removing the maximum element from the fundamental ``double tailed diamond" 
poset  \d$_k(1)$, then it  must be completable to an interval that is isomorphic to all 
of  \d$_k(1)$.
Second, the completing element must cover only elements within that 
interval.
Third, certain kinds of overlaps between two such intervals are prohibited.
It is intriguing that some of the properties obtained have alternate 
derivations in which adding a local structure axiom hypothesis of one of 
these three types allows one of the local structure axiom hypotheses of 
one of the other types to be weakened or omitted.
Further, some of the properties obtained for locally finite posets have alternate derivations in which adding the assumption of finiteness allows one of the local structure axiom hypotheses to be weakened or omitted.
Such trade-offs may parallel interactions among corresponding 
algebraic relations.
More specific comments are made after the statement of Theorem~\ref{k=3 finite theorem}.
We also study the interplay among our axioms in Sections~\ref{k=3 axioms and properties}-\ref{proofs kgeq3} so that we can compare the competing definitions of $d$-complete presented in Section~\ref{d-complete definitions}, and to prepare for proving their equivalence there.

	Much of the structure of a  $d$-complete poset is already 
significantly constrained by the imposition of only the three  $k=3$ 
conditions.
A poset satisfying these is called a  ``$d_3$-complete" poset;  these posets 
may be interesting in their own right. In another sequel to this paper, we visually characterize the 
global and the local structure of finite
$d_3$-complete posets.  This generalizes the classification of  $d$-complete 
posets \cite{DDCT}.

	Section 2 presents the prototypical families of  $d$-complete posets 
and gives additional background, especially for colored $d$-complete posets.
This paper then has three parts with three sections apiece.
The definitions needed for each part are presented in its first 
section and the proofs appear in its last section.
There are notions and results that pertain only to  $d_3$-complete posets,  and the proof 
details for the  $k=3$  conditions are slightly different than the details 
for the  $k\geq4$  conditions.
Therefore Sections~\ref{k=3 axioms and properties},~\ref{k=3 propositions}, and~\ref{proofs k=3} are concerned only with  the axioms needed for $d_3$-complete posets, while Sections~\ref{kgeq3 axioms and properties},~\ref{Propositions kgeq3}, and~\ref{proofs kgeq3} generalize many of those results to  the axioms needed for $d_k$-complete posets for  $k\geq3$.
Section~\ref{d-complete definitions} presents several definitions of  $d$-complete posets and Section~\ref{properties of d-complete posets} presents basic properties of $d$-complete posets. 
Section~\ref{other work} describes appearances of $d$-complete posets. 


\section{Poset terms, prototypical $\bold{\emph{d}}$-complete posets, and more background}

A poset is \emph{locally finite} if every closed interval is finite. `Poset' will mean `locally finite poset' unless `finite' is assumed. Let $P$ be a poset, and let $u,v,w,x,y,z$ denote distinct elements of $P$. We write $x\rightarrow y$ when $y$ covers $x$. We extend this to write $\{x,y\}\rightarrow z$ for $x\rightarrow z$ and $y\rightarrow z$, to write $w\rightarrow\{x,y\}$ for $w\rightarrow x$ and $w\rightarrow y$, and so on. Consult \cite{Enum Comb} for the notions of the (\emph{Hasse}) \emph{diagram} of $P$,  \emph{closed interval} $[w,z]$,  \emph{convex set}, \emph{connected} poset, and \emph{connected components}.
The following definitions come from \cite{CLM}:
A (\emph{covering}) \emph{chain} in $P$ is a set of elements $x_1,\ldots,x_n\in P$ for $n\geq 1$ such that $x_1\rightarrow x_2\rightarrow\cdots\rightarrow x_n$. This chain $C$ has \emph{length} $n-1$; this is denoted $\ell(C)=n-1$. A \emph{rank function} on $P$ is a function $\rho :P\rightarrow \Z$ such that $x\rightarrow y$ implies that  $\rho(y)=\rho(x)+1$. We say $P$ is \emph{ranked} if it has a rank function.

\begin{figure}[b!]
    \centering
   \begin{subfigure}[b]{0.49\textwidth}
        \centering
        \includegraphics[scale=.595]{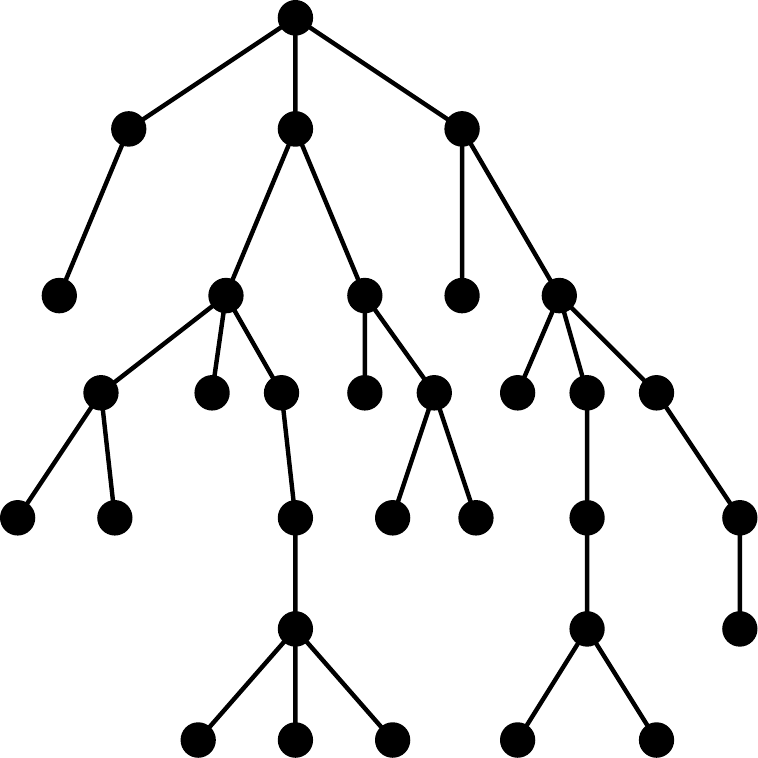} 
            \end{subfigure}
    \hfill
    \begin{subfigure}[b]{0.49\textwidth}
        \centering
        \includegraphics[scale=.2975]{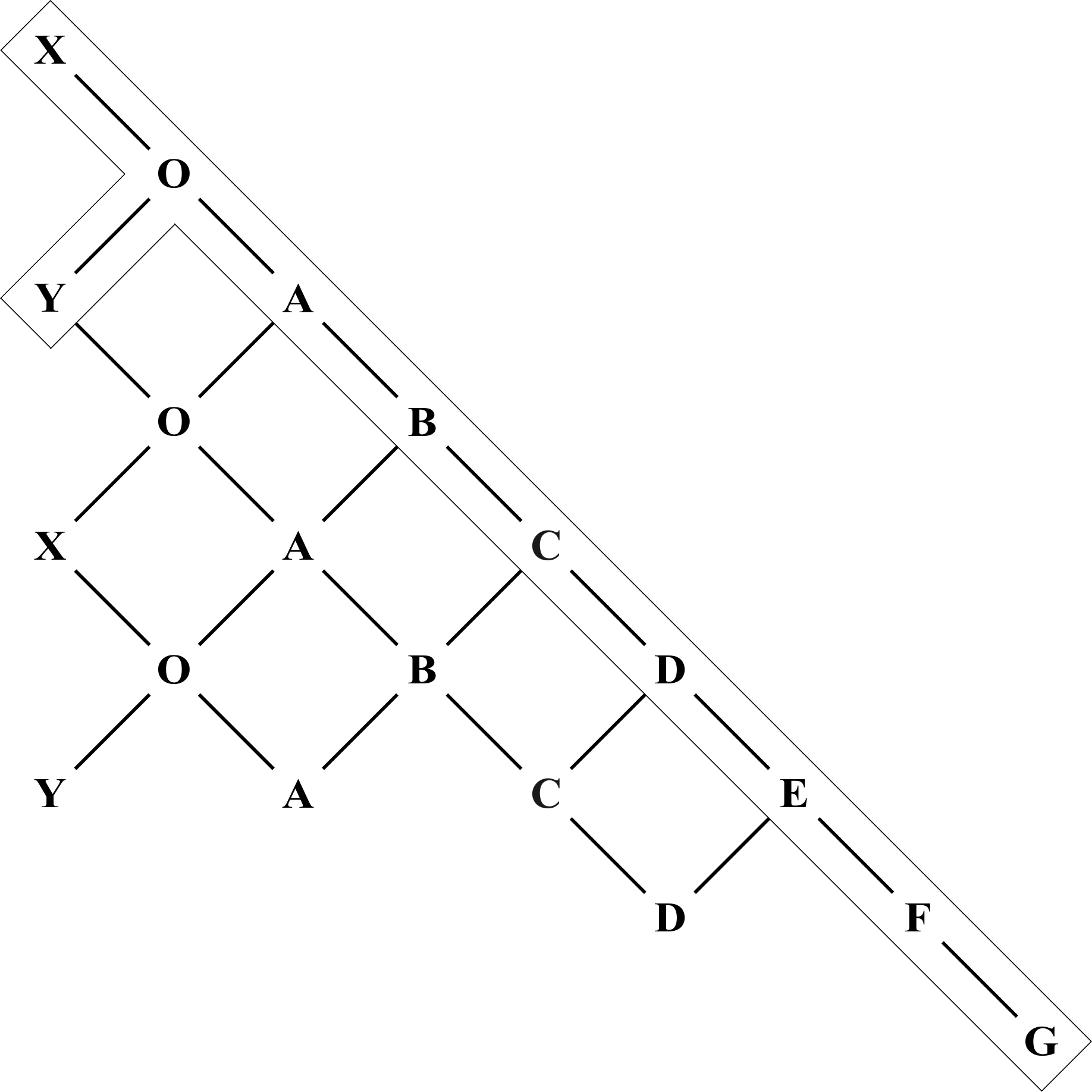} 
    \end{subfigure}
    \vspace{.2in}
    \caption{Rooted tree and shifted shape poset}
   \label{rooted tree shifted shape}
   \end{figure}

\begin{figure}[t!]
    \centering
    \begin{subfigure}[b]{0.49\textwidth}
        \centering
        \includegraphics[scale=.5525]{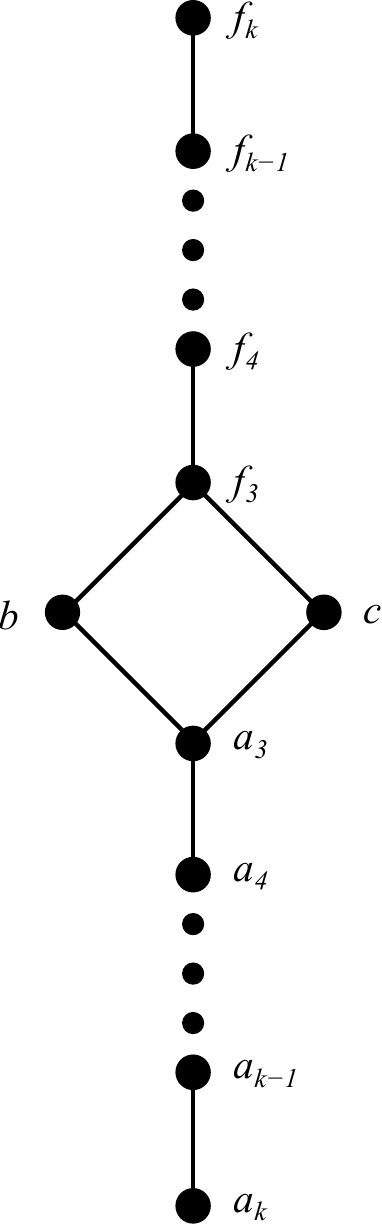} 
            \end{subfigure}
    \hfill
    \begin{subfigure}[b]{0.49\textwidth}
        \centering
        \includegraphics[scale=.5525]{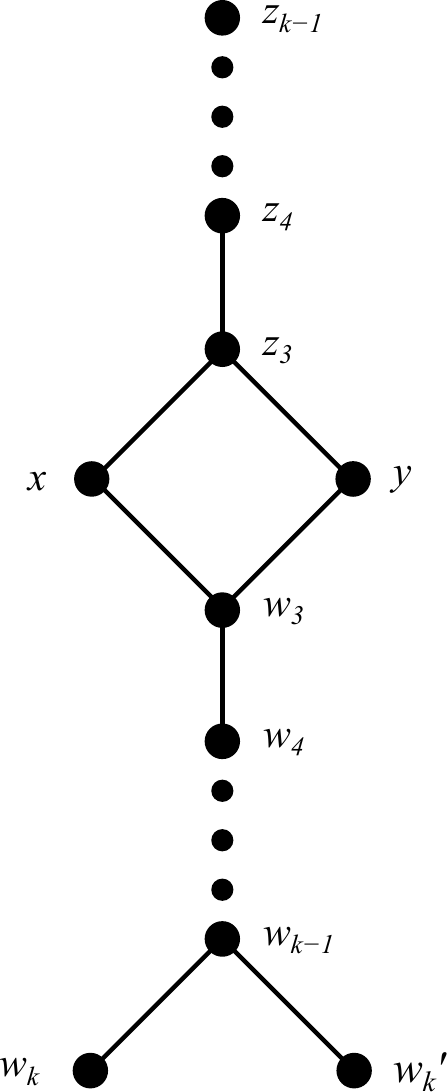} 
    \end{subfigure}
    \vspace{.2in}
    \caption{Double tailed diamond poset \protect\d$_k(1)$ and overlapping $d_k^-$-intervals}\label{dk(1) and overlapping}
   \end{figure}

Definition~\ref{kk} presents the definition of ``$d$-complete" poset. Here we present the prototypical families of posets that motivated the development of that notion:
 A \emph{rooted tree} is a connected poset with a unique maximal element whose diagram is acyclic; see Figure~\ref{rooted tree shifted shape}a.
 Given a Ferrers diagram for a partition of an integer, its \emph{shape} poset is produced by rotating it $45^{\circ}$ clockwise and drawing covering edges between adjacent dots.
\emph{Shifted shape} posets are similarly created from the shifted Ferrers diagrams for strict integer partitions.
Figure~\ref{rooted tree shifted shape}b displays the shifted shape poset for the strict partition $(9,6,3,1)$ of 19. 
Let $k\geq 3$.
The \emph{double tailed diamond} poset \d$_k(1)$ is defined by Figure~\ref{dk(1) and overlapping}a.
Any $(2k-2)$-element self dual poset with exactly two incomparable elements is isomorphic to \d$_k(1)$.
While developing his theory of $P$-partitions, Stanley obtained \cite{OSP} hook length product identities for rooted trees, shapes, and double tailed diamonds that generalized Euler's generating function identity for the number of integer partitions into no more than $n$ parts, and he conjectured such an identity for shifted shapes.

As product identities for bounded $P$-partitions on rectangular shapes and staircase shifted shapes were derived Lie theoretically with Stanley, the Bruhat orders on the quotients $W^J$ of finite Weyl groups that are distributive lattices were classified \cite{BLP}.
Their posets of join irreducible elements  were called ``minuscule" posets; these were labelled with the root system and the dominant weight used to generate the Bruhat order.
Rectangular shapes, staircase shifted shapes, and double-tailed diamonds are minuscule.
Shapes, shifted shapes, the ``filters" of other minuscule posets, and rooted trees are $d$-complete.
Finite connected $d$-complete posets have unique maximal elements.
The \emph{top tree} of a finite connected poset with a unique maximal element is the rooted tree that consists of the elements $x$ such that $\{y\colon y\geq x\}$ is a chain.
The shifted shape in Figure~\ref{rooted tree shifted shape}b is a filter of the staircase minuscule poset denoted \d$_{10}(\omega_9)$; its top tree is circled.
The top tree of a filter of a minuscule poset is the Dynkin diagram for the Weyl group from which its Bruhat order was formed.
In the shifted shape the top tree is the Dynkin diagram \D$_{10}$.
For $n\geq 3$ the double tailed diamond \d$_n(1)$ is the minuscule poset \d$_n(\omega_1)$, and its top tree is the Dynkin diagram \D$_n$. Hence for $k\geq 3$ the subscript `$k$' in \d$_k(1)$ is the number of generators for the associated Weyl (Coxeter) group of type \D$_k$.

\begin{figure}[h!]
       \centering
    \begin{subfigure}[c]{0.49\textwidth}
        \centering
                \includegraphics[scale=.2975]{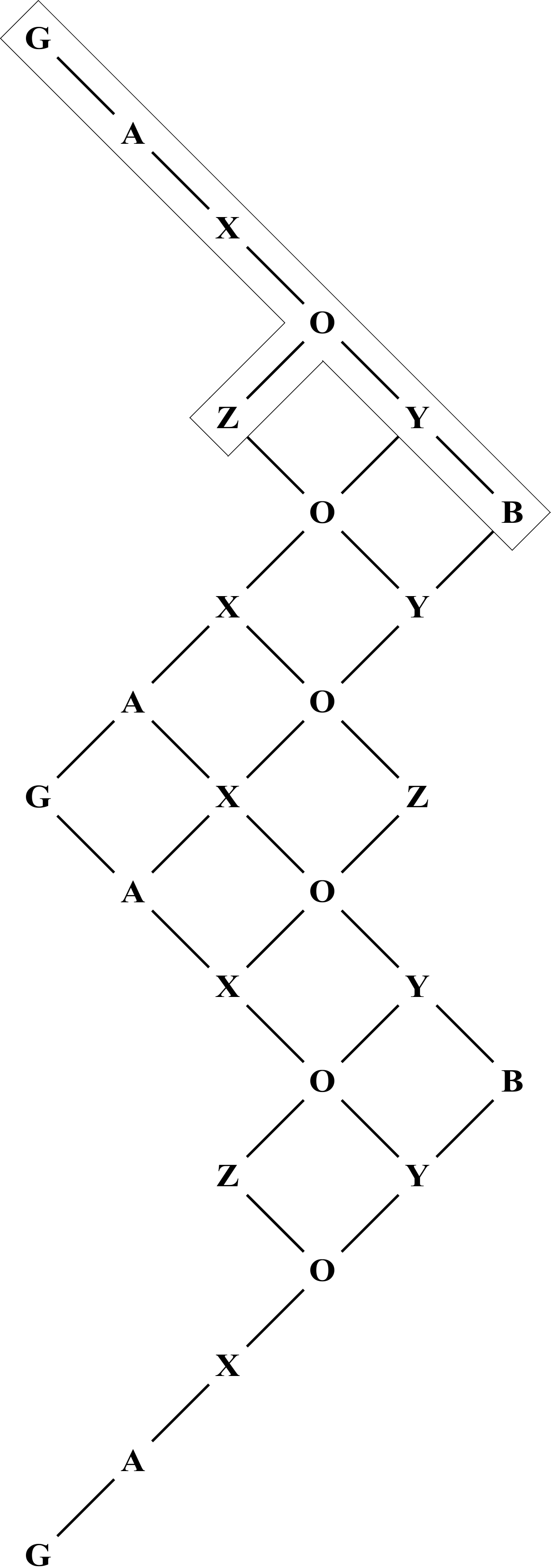} 
            \end{subfigure}
    \hfill
    \begin{subfigure}[c]{0.49\textwidth}
        \centering
        \includegraphics[scale=.255]{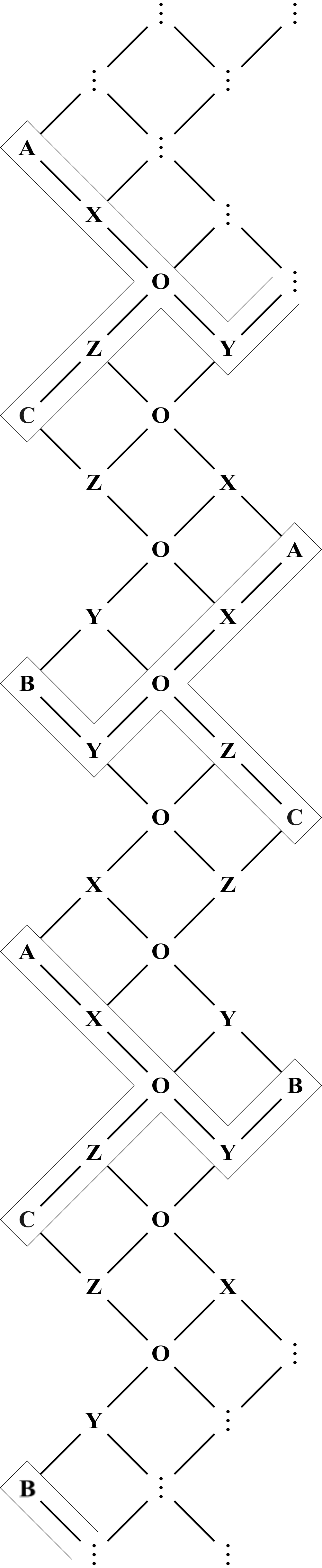} 
    \end{subfigure}
    \vspace{.2in}
    \caption{Minuscule poset \protect\e$_7(7)$ and the poset for $\widetilde{\text{\protect\E}}_6(\omega_A-\omega_B)$}
    \label{e7(7) E6}
   \end{figure}

A finite connected poset with a unique maximal element is \emph{simply colored} \cite{Wave} if its elements have been colored such that the elements in its top tree are distinctly colored, every other element receives one of those colors, and no two elements in a chain interval receive the same color.
Inspired by a formulation \cite{Ste} of Stembridge, we now \cite{Japan} define it to be a \emph{(simply) colored $d$-complete} poset if: equichromatic elements are comparable, any two elements with colors that are adjacent in the top tree are comparable, the colors of two elements that are adjacent in the Hasse diagram are adjacent in the top tree, and in the open interval between two consecutive equichromatic elements there are exactly two elements whose colors are adjacent to that color in the top tree.
Figures~\ref{rooted tree shifted shape}b and~\ref{e7(7) E6}a display such colored posets; the latter one is the exceptional minuscule colored $d$-complete poset \e$_7(7)$.
The notion of colored $d$-complete poset was developed \cite{Mar} when a combinatorial linear algebra version of the geometric representation basis theorem of Seshadri that had been used \cite{BLP} to prove $m$-bounded $P$-partition identities  was developed.
For finite posets, the notions of $d$-complete poset and of colored $d$-complete poset are essentially equivalent \cite{Wave}: Ignoring the colors of one of the latter produces one of the former. And given one of the former, its elements may be colored in essentially only one way to produce one of the latter. (The (colored) $d$-complete definition used in \cite{Wave} was the order dual of the definition used here and elsewhere.)



Dale Peterson introduced \cite{Carrell} the notion of a ``$\lambda$-minuscule"element $w$ of a Kac-Moody Weyl group $W$.
When $W$ is simply laced, it was shown \cite{Wave} for such an element that the ``ideal" $(w)$ of the Bruhat order on $W$ is a distributive lattice.
It was further shown that the order dual of a poset $P$ is colored $d$-complete if and only if $P$ arises as the poset of join irreducible elements of such a distributive lattice $(w)$ for some dominant $\lambda$. 
Here $(w)\cong J(P)$ in the language of \cite{Enum Comb}.
Then, as the ``heap" of $w$, the colored poset contains much information \cite{Ste} concerning the reduced decompositions of $w$.

\begin{figure}[h!]
    \centering
    \begin{subfigure}[c]{0.49\textwidth}
        \centering
        \includegraphics[scale=.425]{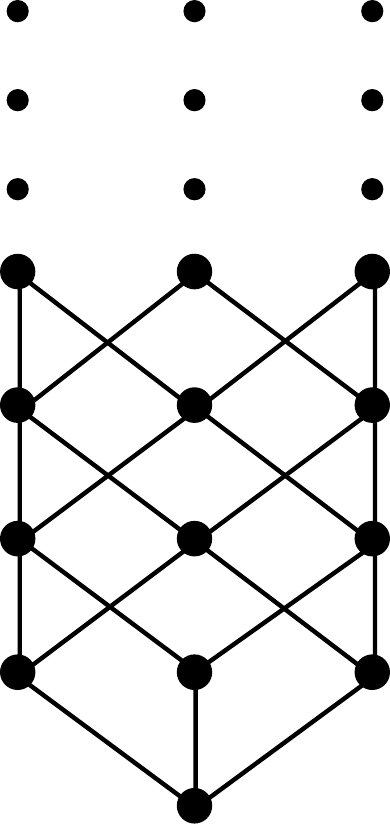} 
            \end{subfigure}
    \hspace{-.8in}
    \begin{subfigure}[c]{0.49\textwidth}
        \centering
        \includegraphics[scale=.425]{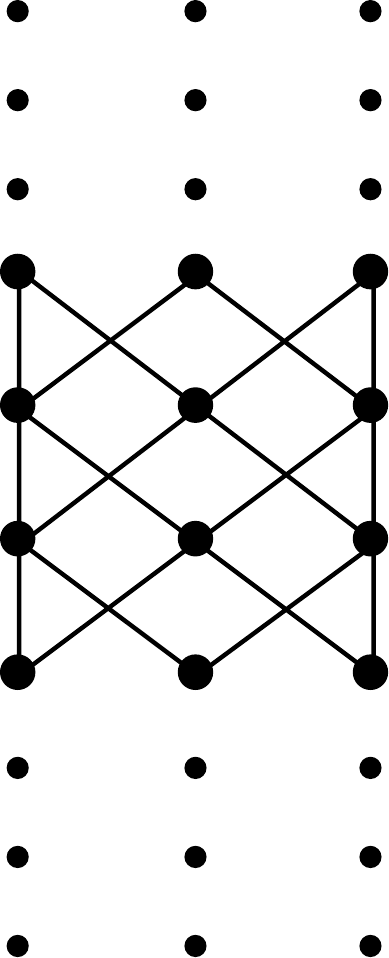} 
    \end{subfigure}
    \vspace{.2in}
    \caption{Locally finite posets}
    \label{infinite d-complete posets}
   \end{figure}

The definitions of ``$d$-complete poset" require that the intervals in the poset that are isomorphic to \d$_k(1)$ (or are nearly isomorphic to \d$_k(1)$) for $k\geq 3$ are well behaved in certain respects.
So one may view a $d$-complete poset as consisting of double tailed diamonds that have been carefully ``woven" together.
The doubly infinite colored poset displayed in Figure~\ref{e7(7) E6}b is generated from the weight $\lambda=\omega_{A}-\omega_{B}$ for the exceptional affine Weyl group $\widetilde{\text{\E}}_6$ in the ``numbers game" manner implicitly used in the proofs of Lemma 3.2 and Proposition 3.1 in \cite{BLP}.
In addition to satisfying our definition of $d$-complete for locally finite uncolored posets, this colored poset also satisfies the requirements above for finite simply colored $d$-complete posets, once the role of top tree has been taken by an embedded copy of the Dynkin diagram.
The intervals formed from the consecutive occurrences of a color in Figure~\ref{rooted tree shifted shape}b and Figure~\ref{e7(7) E6} are double tailed diamonds.
In this paper both of the infinite uncolored posets in Figure~\ref{infinite d-complete posets} qualify to be $d$-complete posets; see the Added Note.

Roughly speaking, the overall global structure of a finite connected $d$-complete poset is that of a rooted tree, but with interspersed ``slant irreducible" components \cite{DDCT}.
These irreducible components fall into 15 classes, of which 14 are indexed by top trees which are Dynkin diagrams of general type \E.
The sole member of the 15th class is \e$_7(7)$.



\section{Definitions, axioms, and properties for $\bold{\emph{k}=3}$}
\label{k=3 axioms and properties}

Let $P$ be a poset with distinct elements $u,w,w',x,y,z$.  A \emph{diamond} is a subset $\{w;x,y;z\}$ of $P$   such that  $w\rightarrow \{x,y\}$ and $\{x,y\}\rightarrow z$. The \emph{bottom} and \emph{top} of the diamond are $w$ and $z$ respectively, and the \emph{elbows} are $x$ and $y$.  
We say the top z is \emph{free} if it covers only $x$ and $y$.
An interval $[w,z]$ is a \emph{$d_3$-interval} if it is a diamond $ \{w;x,y;z\}$ for some $x,y\in P$. 
Once it is known that a diamond $\{w;x,y;z\}$ forms all of the interval $[w,z]$, we refer to $z$ as the maximum element.
A subset $\{w;x,y\}$ of $P$ is  a \emph{vee} or a $d_3^-$-\emph{set} if  $w\rightarrow \{x,y\}$. Note that a $d_3^-$-set is convex.
A $d_3^-$-set $\{w;x,y\}$ is \emph{completed} if there exists a \emph{completing element} $z$ such that $\{w;x,y;z\}$ is a $d_3$-interval. 
Two $d_3^-$-sets $\{w;x,y\}$ and $\{w';x,y\}$ are said to \emph{overlap}.
Two $d_3$-intervals (or diamonds) are \emph{distinct} if their set symmetric difference is non-empty.
An interval $[w,z]$ is \emph{short} if there exists $u$ such that $w\rightarrow u\rightarrow z$.

\begin{fact}\label{get interval}
Let $S=\{w;x,y\}$ be a $d_3^-$-set and let $z\in P$. 

\noindent (a) If $S\cup\{z\}$ is a $d_3$-interval, then $S\cup\{z\}=[w,z]$ and  $z$ covers $x$ and $y$. 

\noindent (b) If $z$ covers $x$ and $y$ and no other elements, then $S\cup\{z\}$ is the $d_3$-interval $[w,z]$.
\end{fact}

We call a structural property an ``axiom" if  we use it as part of a definition of the $d$-complete property in Section~\ref{d-complete definitions} or near the end of Section~\ref{k=3 propositions}.

\noindent \textbf{Axioms}  \ \textit{A poset satisfies the}

\vspace{.1in}

\noindent [\textit{Class I: Completion Axioms}]

\begin{axiom*}\hspace{-.075in}(VT)  \hspace{.001in}
``Vees have Tops" axiom if  $w\rightarrow \{x,y\}$ implies that there exists $z$ such that $\{x,y\}\rightarrow z$,  the
\end{axiom*}

\begin{axiom*}\hspace{-.075in}($D3^-C$)  \hspace{.001in}
``$d_3^-$-sets are Completed" axiom if for each $d_3^-$-set $S$ there exists $z$ such that $S\cup \{z\}$ is a $d_3$-interval, the
\end{axiom*}

\vspace{.1in}

 \noindent [\textit{Class II: Freeness Axioms}]

\begin{axiom*}\hspace{-.075in}(FT)  \hspace{.001in}
 ``diamonds have Free Tops" axiom if the top element of each diamond covers only the elbows of the diamond, the
 \end{axiom*}
 
 \begin{axiom*}\hspace{-.075in}(D3MF)  \hspace{.001in}
``$d_3$-interval Maxs are Free" axiom if the maximum element of each $d_3$-interval covers only the elbows of the interval, the
\end{axiom*}

\vspace{.1in}

 \noindent[\textit{Classes I/II: Completion/Freeness Axiom}]

\begin{axiom*}\hspace{-.075in}($D3^-CF$)  \hspace{.001in}
 ``$d_3^-$-sets are Completed Freely" axiom if for each $d_3^-$-set $S$ there exists $z$ such that $S\cup \{z\}$ is a $d_3$-interval with maximum element $z$ such that $z$ covers only elements from $S$, the
 \end{axiom*}

\vspace{.1in}

\noindent[\textit{Class III: Forbidden Structure Axioms}]

\begin{axiom*}\hspace{-.075in}(NCC)  \hspace{.001in}
 ``No Criss Cross" axiom if there do not exist overlapping $d_3^-$-sets, and the
 \end{axiom*}

\begin{axiom*}\hspace{-.075in}(D3MD)  \hspace{.001in}
``$d_3$-interval Maxs are Distinct" axiom if the maximum elements of distinct $d_3$-intervals are distinct.
\end{axiom*}

\noindent At times, Axioms VT and FT together are referred to as Diamond I+II and Axioms D3$^-$C and D3MF together are referred to  as Classic I+II.

\begin{myremark}\label{obvious axiom implications k=3} 
Axioms D3$^-$C, D3$^-$CF, and FT obviously respectively imply VT, D3$^-$C, and D3MF. Also, Axioms D3$^-$C and D3MF together imply D3$^-$CF. 
\end{myremark}

 \noindent\textbf{Properties}. \textit{A poset has the}

 \begin{property*}\hspace{-.075in}(UPUE) \hspace{.001in}
 ``Upward Propagation of Up Edges" property if  $w\leq y$, $w\rightarrow x$, and $x\not\leq y$ imply there exists $z$ such that $y\rightarrow z$ and $x\leq z$, the
 \end{property*}

\begin{property*}\hspace{-.075in}(UM) \hspace{.001in}
``Unique Maximal element" property if it has a unique maximal element, the
\end{property*}

\begin{property*}\hspace{-.075in}(CLMEE) \hspace{.001in}
``Chain Lengths to Maximal Element are Equal" property if it has UM and if every chain from an element $w$ to the unique maximal element $z$ has the same length, the
\end{property*}

\begin{property*}\hspace{-.075in}(CLE)   \hspace{.001in}
``Chain Lengths are Equal" property if whenever $x<y$ all chains from $x$ to $y$ have equal length, the
\end{property*}

\begin{property*}\hspace{-.075in}(SS)  \hspace{.001in}
``Short Intervals are Small" property if whenever $[w,z]$ is a short interval, then $|[w,z]|\in\{3,4\}$, the
\end{property*}

\begin{property*}\hspace{-.075in}(DAI)  \hspace{.001in}
``Diamonds Are Intervals" property if each diamond is an interval, the
\end{property*}

\begin{property*}\hspace{-.075in}(NTC)  \hspace{.001in}
``No Triply Covereds" property if no element is covered by three elements, the
\end{property*}

\begin{property*}\hspace{-.075in}(UT)  \hspace{.001in}
 ``Unique Top" property if each vee has exactly one top, and the
 \end{property*}

\begin{property*}\hspace{-.075in}(UC3)  \hspace{.001in}
``Unique Completion" property if each $d_3^-$-set has exactly one completing element.
\end{property*}

\noindent The ranked and connected properties are indicated with the labels (Rank) and (Conn).

\begin{myremark}\label{NCC axiom}
If Property NTC were to be regarded as an axiom, it would belong to Class III since it could be viewed as prohibiting two more kinds of overlap between two $d_3^-$-sets that are not prohibited by Axiom NCC:
Let $\{w;x,y\}$ and $\{w';x',y'\}$ be two distinct $d_3^-$-sets. Suppose these sets have a coincidence between their minimal elements and/or a coincidence among their maximal elements. If $\{x,y\}=\{x',y'\}$, then $w\neq w'$ and this is prohibited by NCC. If $\{x,y\}\neq\{x',y'\}$ and $w=w'$, then Property NTC prohibits $w=w'$ from being covered by three or four distinct elements from $\{x,y,x',y'\}$. (Axiom NCC and Property NTC together do not prohibit the one remaining possibility, the ``W": here $|\{x,y\}\cap\{x',y'\}|=1$ and $w\neq w'$.)
\end{myremark}

It is easy to see that:
\begin{fact}[\normalfont{DAI $\#1$}] \hspace{.001in} \label{DAI from NTC} If a poset $P$ is FT or has No Triply Covereds or has Short Intervals are Small, then it has Diamonds are Intervals.
\end{fact}


\section{Results for $\bold{\emph{k}=3}$}
\label{k=3 propositions}

\vspace{-.08in}
 The implications in Theorem~\ref{k=3 locally finite theorem} below are displayed in tabular form, with the first four columns displaying hypotheses that are Class I, II, III axioms or a property.
 For example, Parts (\classiccombo), (\MaxsDistinctMaxsFree), and (\secondNCCMaxsFree) have the acronym D3$^-$CF entered midway  between the columns for Class I and Class II axioms.
 This indicates that the Class I/II hybrid axiom ``$d_3^-$-sets are Completed Freely" is being assumed in these parts.
 Parts (\MaxsDistinctMaxsFree)  and (\secondNCCMaxsFree) together indicate that the Class III axioms ``No Criss Cross" and ``$d_3$-intervals Maxs are Distinct" are equivalent in the presence of D3$^-$CF.
 Parts (\firstcompleted)-(\diamondtoclassic) describe ways in which the easy-to-check ``Vees have Tops" axiom may be strengthened to the ``$d_3^-$-sets are Completed" axiom needed for a $d$-complete poset.
 Five of the parts are concerned only with axioms; this illustrates the interplay among the axioms mentioned in the introduction.
 For example, part of Part (\MaxsDistinctMaxsFree) strengthens the weak Class II requirement contained in the Class I/II hybrid axiom D3$^-$CF  to the full-strength Class II axiom D3MF when the Class III axiom NCC is present.
 And without a Class I axiom being present, in Part (\MaxsDistinct) the Class III axiom NCC implies the other Class III axiom D3MD when the Class II axiom D3MF is present.
An entry of ``etc." in the last column indicates that further conclusions may be drawn using one or both of the listed conclusions to satisfy an earlier line in the table.
 The last part propagates an edge in a vee upwardly along a chain when VT is present.


\begin{theorem} \label{k=3 locally finite theorem} The implications in Table~\ref{k=3 locally finite theorem table} hold in a poset. \\
\begin{center}
\vspace{-.102in}
\resizebox{1.02\totalheight}{!}{
\begin{tabular}{ c  @{\hspace{.4in}}c   @{\hspace{-.1in}}c   @{\hspace{-.1in}}c  c c  c c  c c}

&  \underline{I}  & & \underline{II} &  &\underline{III} & & \underline{Property} & &\underline{Conclusion(s)} \\

(\UT)&  VT && --- && NCC && --- && UT \\	
(\firstcompleted)   & VT & &--- & &--- & & DAI & & $D3^-C$ \\
(\secondcompleted)   &  VT && ---&& ---& & NTC && $D3^-C$ \\
(\thirdcompleted)  &   VT && ---&& ---&& SS && $D3^-C$ \\
(\diamondtoclassic) & VT && FT && ---&& ---&& $D3^-C$ + D3MF \\
(\secondDAI) &  $D3^-C$ && ---&& ---&& UT && DAI \\
(\DAIUC) &  $D3^-C$ && ---&& NCC &&---&& UT, DAI, UC3 \\
(\NCCUC)  &  $D3^-C$ &&---&& D3MD && NTC && NCC, etc.\\
(\MaxsDistinct) & ---&&  D3MF && NCC && ---&& D3MD \\
(\FT) & ---&&  D3MF &&--- && DAI && FT \\
(\classiccombo)  && $D3^-CF$ & &&--- && UC3 && $D3^-C$ + D3MF \\
(\MaxsDistinctMaxsFree)  & & $D3^-CF$ && & NCC &&---& & D3MF, D3MD, etc. \\
(\secondNCCMaxsFree)   & & $D3^-CF$ && & D3MD && ---&& NCC, D3MF, etc. \\
(\fullclassictodiamond) &  $D3^-C$ && D3MF && NCC &&--- && VT + FT, etc.\\
(\UPUE) & VT &&--- && ---&& ---&&  UPUE \\
\end{tabular}}
\captionof{table}{Implications for Theorem~\ref{k=3 locally finite theorem}}
\label{k=3 locally finite theorem table}
\end{center}
\end{theorem}

\noindent In Section 9 we see that the hypotheses of Parts (\MaxsDistinctMaxsFree), (\secondNCCMaxsFree), and (\fullclassictodiamond) satisfy the definition of ``$d_3$-complete" poset.
Posets satisfying these axioms satisfy all of the $k=3$ axioms and have the DAI, UT, and UC3 properties.
Given this remark, it can be seen that Part (\fullclassictodiamond) is closely related to Part (\MaxsDistinctMaxsFree).
We have included Part (\fullclassictodiamond) because it gives a partial converse to Part (\diamondtoclassic), and because it clarifies the misworded statement ``We have just required~$\ldots$" on pp. 65 and 283 of \cite{DDCT} and \cite{Wave}; that statement should have instead begun ``It can be shown that~$\ldots$".

Combining Remark~\ref{obvious axiom implications k=3}, Fact~\ref{DAI from NTC}, and Part (\FT) of Theorem~\ref{k=3 locally finite theorem}, we have:

\begin{mycorollary}
A poset is FT if and only if it is D3MF and has Diamonds are Intervals.
\end{mycorollary}

Now we consider finite posets. We believe that the converse of Part (\thirdcompleted) above holds here:

\begin{conj}\label{ss conjecture} If a finite poset is $D3^-C$, then it has Short intervals are Small (and is VT).
\end{conj}

\noindent The next theorem presents some  implications that may be deduced when the poset is finite. 
Parts (\UMCLMEE) and (\RankCLE) present four fundamental structural properties that follow from the Upward Propagation of Up Edges part above. 
Part (\classictodiamond) says that the converse of Theorem~\ref{k=3 locally finite theorem}(\diamondtoclassic)  is known to hold when the poset is finite.
The first part of Part (\thirdNCCsecondNTCUC) says that the D3$^-$CF hypothesis for the first part of Theorem~\ref{k=3 locally finite theorem}(\secondNCCMaxsFree)  may be weakened to D3$^-$C in the finite case.

 \begin{theorem} \label{k=3 finite theorem} The implications in Table~\ref{k=3 finite theorem table} hold in a finite poset. \\
 \begin{center}

\begin{tabular}{ c  @{\hspace{.4in}}c   @{\hspace{-.1in}}c   @{\hspace{-.1in}}c  c c  c c  c c}

&  \underline{I}  & & \underline{II} &  &\underline{III} & & \underline{Property} & &\underline{Conclusion(s)} \\
(\UMCLMEE)   &  VT & &---  & & --- & & Conn & & UM, CLMEE\\
(\RankCLE)  &  VT && --- && --- &&  --- && Ranked, CLE \\
(\Ss) & VT&& --- &&  --- && NTC && SS \\
(\firstNTC) & & $D3^-CF$ && & ---  && --- && NTC \\
(\classictodiamond) &  $D3^-C$ && D3MF && --- &&---  && VT + FT \\
(\thirdNCCsecondNTCUC) &  $D3^-C$ &&  --- && D3MD && --- && NCC, NTC, UC3, etc.\\
\end{tabular}
\captionof{table}{Implications for Theorem~\ref{k=3 finite theorem}}
\label{k=3 finite theorem table}
\end{center}
\end{theorem}

\begin{figure}[t!]
    \centering
    \begin{subfigure}[c]{0.49\textwidth}
        \centering
        \includegraphics[scale=.6]{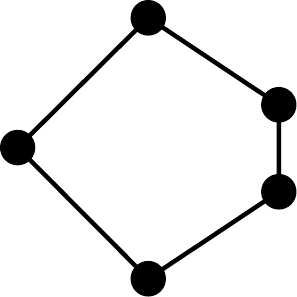}
            \end{subfigure}
    \hspace{-.8in}
        \begin{subfigure}[c]{0.49\textwidth}
        \centering
        \includegraphics[scale=.6]{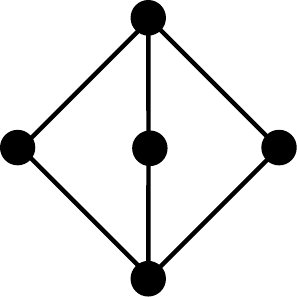}
    \end{subfigure}
    \vspace{.2in}
    \caption{Seeds for proposed counterexamples}
    \label{seeds}
   \end{figure}

Figure~\ref{infinite d-complete posets}a gives a counterexample to dropping the assumption of finiteness from Part (\firstNTC) and from the second part of Part (\thirdNCCsecondNTCUC) of this theorem.
 Requiring Axioms D3$^-$C or D3$^-$CF without requiring the completions of vees to be unique can lead to messy situations if insufficient requirements have been imposed with Class II or Class III axioms or with finiteness.
 But we have not considered Property UC3 as an axiom in this paper since uniqueness can be difficult to confirm.
 We believe that Parts (\Ss) and (\classictodiamond) and the first part of Part (\thirdNCCsecondNTCUC) of this theorem also do not hold when finiteness is dropped.
 Figure~\ref{seeds}a presents a ``seed" for a proposed counterexample to Part (\Ss) and Figure~\ref{seeds}b does so for proposed counterexamples to Part (\classictodiamond) and  the first part of Part (\thirdNCCsecondNTCUC).
 The failure of the diagrams generated from these seeds to ``close up" in a neat finite fashion may correspond to some kind of messy infinite algebraic quotient that has resulted from insufficient relations having been imposed.
 In a similar vein, attempting to prove that the contrapositive ``$\neg$ SS $\Rightarrow \neg$ D3$^-$C $\lor \neg$ finite" of Conjecture~\ref{ss conjecture} is true when the second diagram of  Figure~\ref{seeds}b is present as a subdiagram also seems to generate an infinite poset.

\begin{conj}
There exist infinite locally finite posets that are counterexamples to the parts of Theorem~\ref{k=3 finite theorem} mentioned above.
\end{conj}

\begin{myremark} Since VT and FT imply D3$^-$C and D3MF, and D3$^-$C and D3MF obviously imply D3$^-$CF, Theorem~\ref{k=3 finite theorem}(\firstNTC)  implies that a finite poset has NTC whenever it is Diamond I+II or Classic I+II.
\end{myremark}

Several definitions of $d_3$-complete will be given for locally finite posets in Section~\ref{d-complete definitions}.
For this paragraph, let us use that Classic definition to say that a poset is $d_3$-complete if it is D3$^-$C, D3MF, and NCC.
This provides a context to discuss the interplay among the axioms and between the axioms and the assumption of finiteness, especially in regard to forming other combinations of axioms that are equivalent to the Classic definition.
Within Class I, Axiom D3$^-$C is stronger than VT.
Within Class II, Axiom FT is stronger than D3MF.
It is interesting that in the diamond point of view, using the stronger FT compensates for using the weaker VT in Theorem~\ref{k=3 locally finite theorem}(\diamondtoclassic) so that one can still obtain the combination D3$^-$C plus D3MF needed in the $d_3$-interval point of view for the Classic definition.
Conversely, in the $d_3$-interval point of view when NCC is present (Theorem~\ref{k=3 locally finite theorem}(\fullclassictodiamond)) or the poset is finite (Theorem~\ref{k=3 finite theorem}(\classictodiamond)), using the stronger D3$^-$C compensates for using the weaker D3MF so that one can still obtain the more convenient  combination of VT plus FT in the diamond point of view. 
What happens if the weaker Class~I axiom VT is paired with the weaker Class~II axiom D3MF?
By Theorem~\ref{k=3 locally finite theorem}(\secondcompleted), strengthening the Class III axiom of NCC by also assuming NTC to prohibit two more same-rank overlaps between two $d_3^-$-sets allows one to satisfy the Classic definition of $d_3$-complete with the combination of VT, D3MF, NCC, and NTC. (This was done in some earlier editions of \cite{JDT}.)
When the poset is finite, Remark~\ref{obvious axiom implications k=3} and Theorem~\ref{k=3 finite theorem}(\firstNTC)  indicate that this strengthening did not go too far: here a finite poset that is $d_3$-complete by the Classic definition already has the NTC property.



\section{Proofs for $\bold{\emph{k}=3}$}
\label{proofs k=3}

\begin{proof} \hspace{-.035in}\textit{of Theorem~\ref{k=3 locally finite theorem}}.
Parts (\UT),(\firstcompleted), (\FT), and (\classiccombo) follow quickly from the definitions; Part (\MaxsDistinct) also follows directly with a bit of thought. 
For Parts (\secondcompleted) and (\thirdcompleted), note that: If VT is present, then D3$^-$C can fail at a $d_3^-$-set only if there is an ``extra" chain from the minimum element to the diamond top required by VT. Such a chain is ruled out by either NTC or SS. 
For Part (\diamondtoclassic), recall from Section~\ref{k=3 axioms and properties} the implications FT $\Rightarrow$ D3MF and FT $\Rightarrow$ DAI. Then by Part (\firstcompleted) we get D3$^-$C.

For Part (\secondDAI), given a diamond $\{w;x,y;z\}$, for $\{w;x,y\}$ by D3$^-$C there is a $z'$ such that $[w,z']$ is a $d_3$-interval. So UT implies $z'=z$, implying that the diamond is an interval. For Part (\DAIUC), follow Remark~\ref{obvious axiom implications k=3} and Part (\UT) by Part (\secondDAI).

For Part (\NCCUC), suppose $\{w,w'\}\rightarrow\{x,y\}$. Applying D3$^-$C to $\{w;x,y\}$ gives $z$ such that $[w,z]$ is a $d_3$-interval. Here D3MD implies that $[w',z]$ is not a $d_3$-interval. So there exists some $u\in[w',z]$ such that $u\not\in\{w';x,y;z\}$. Let $v$ be the minimal such element; we have $w'\rightarrow v$. Since $v\not\in\{x,y\}$, this would violate NTC. So NCC holds. 

For Part (\MaxsDistinctMaxsFree), first note D3$^-$CF $\Rightarrow$ D3$^-$C; add NCC with Part (\DAIUC) to obtain UC3, which via Part (\classiccombo) implies D3MF. Then Part (\MaxsDistinct) gives D3MD.
 For Part (\secondNCCMaxsFree), suppose $\{w,w'\}\rightarrow\{x,y\}$. Applying D3$^-$CF to $\{w;x,y\}$ gives $z$ such that $[w,z]$ is a $d_3$-interval with $z$ free. This $z$ covers exactly $x$ and $y$ in $\{w';x,y\}$. So $[w',z]$ is a $d_3$-interval. This contradicts D3MD at $z$.
For Part (\fullclassictodiamond), first note D3$^-$C $\Rightarrow$ VT; add NCC with Part (\UT), which via Part (\secondDAI) gives DAI. Then Part (\FT) gives FT.

 Part (\UPUE) was Proposition~F1 of \cite{DDCT}, which did not actually need finiteness.
\qed
\end{proof}

\begin{proof}\hspace{-.035in}\textit{of Theorem~\ref{k=3 finite theorem}}.
 Part (\UMCLMEE) was Propositions~F2 and~F3 of \cite{DDCT}.
For Part (b), each component has a unique maximal element. Then CLMEE can be used to construct a well-defined rank function on each component. Property CLE follows.

For Part (\Ss), a short interval $[w,z]$ will be contained in a connected component. Let $u$ be such that $w\rightarrow u\rightarrow z$. Suppose $|[w,z]|\geq5$. So there exists $x,y\in[w,z]$ such that $w,u,x,y,z$ are distinct. There exist chains from $w$ to $z$ that pass through $x$ and through $y$. Here CLE implies that these chains are of length 2. So $w\rightarrow\{u,x,y\}$, contradicting NTC.

For Part (\firstNTC), suppose $w\rightarrow\{x_1,y_1,z_1\}$. Applying D3$^-$CF three times yields three distinct free completing elements $x_2, y_2, z_2$. This axiom can be repeatedly applied three times in this fashion ad infinitum, contradicting finiteness.
For Part (\classictodiamond), note that D3$^-$C implies VT, and adding in D3MF gives D3$^-$CF. Part (\firstNTC) provides NTC, which was used in Fact~\ref{DAI from NTC} to get DAI. Then Theorem~\ref{k=3 locally finite theorem}(\FT)  gives FT.

To prove the first part of  Part (\thirdNCCsecondNTCUC), suppose $\{x_0,y_0\}\rightarrow\{x_1,y_1\}$. Applying axioms D3$^-$C and D3MD together twice implies that there exist distinct completing elements $x_2$ and $y_2$. These axioms can be repeatedly applied twice in this fashion ad infinitum, contradicting finiteness. The proof of the second part of Part (\thirdNCCsecondNTCUC) is the same as the proof of Part (\firstNTC), except now D3$^-$C and D3MD are used instead of D3$^-$CF to construct the infinite succession of completing elements three at a time. For the third part use Theorem~\ref{k=3 locally finite theorem}(\DAIUC).
\qed
\end{proof}

\vspace{.1in}
\section{Definitions, axioms, and properties for $\bold{\emph{k}\geq3}$}
\label{kgeq3 axioms and properties}
\vspace{.02in}

 Let $P$ be a poset.  Here and below $u,v,w,x,y$ denote arbitrary elements. We define convex sets $[w;x,y]\defeq[w,x]\cup[w,y]$ and $[u,v;x,y]\defeq[u;x,y]\cup[v;x,y]$.
 \vspace{.1in}

Let $k\geq 3$. Consider the double tailed diamond (DTD) poset \d$_k(1)$ of Figure~\ref{dk(1) and overlapping}a. The two incomparable elements $b$ and $c$ are its \emph{elbows}.  Its \emph{neck elements} are $f_3,f_4,\ldots,f_k$ and its \emph{tail elements} are $a_3,a_4,\ldots,a_k$. When $k\geq 4$, all but the lowest of the neck elements are its  \emph{strict} neck elements and all but the highest of the tail elements are its \emph{strict} tail elements. 
A $Y_k$\textit{-set} $[w;x,y]\subseteq P$ is a convex set such that $[w;x,y]\cong [a_k;b,c]\subseteq \text{\d}_k(1)$. Note that a $Y_3$-set is a vee.
Suppose the elements of this set are $w\eqdef w_k\rightarrow w_{k-1}\rightarrow\cdots\rightarrow w_3\rightarrow\{x,y\}$.
Here $w_k,\ldots,w_3$ are the \textit{stem} elements of $Y_k$.
A $\Lambda Y_k$\textit{-set} $[u,v;x,y]\subseteq P$ is a convex set of the following form: we require $\{u,v\}\rightarrow w_k$ and $[u,v;x,y]=\{u,v\}\cup[w_k;x,y]$, where $[w_k;x,y]$ is a $Y_k$-set.
 
\vspace{.1in}
An interval $[w,z]$ in $P$ is a \emph{$d_k$-interval} if $[w,z]\cong \text{\d}_k(1)=\{a_k,\ldots,a_3;b,c;f_3,\ldots,f_k\}$. We say that $[w,z]$ is a \emph{DTD} interval if we do not want to mention k. Note that if $[w,z]=\{w\eqdef w_k,\ldots,w_3;x,y;z_3,\ldots z_k\defeq z\}$ is a $d_k$-interval, then $[w_h,z_h]$ is a $d_h$-interval for $3\leq h\leq k$. 
A neck element $u$ in a $d_k$-interval $[w,z]$ is \emph{free} if $u$ covers only (an) element(s) in $[w,z]$.
A convex set is a $d_k^-$\textit{-set} if it is isomorphic to \d$_k(1)-\{f_k\}$.
When $k\geq 4$, a $d_k^-$-set is an interval and thus may be referred to as a $d_k^-$\emph{-interval}.
Here a $d_k^-$-interval $[w,z']$ may be described as $\{w\eqdef w_k,w_{k-1},\ldots,w_3;x,y;z_3,\ldots,z_{k-2},z_{k-1}\defeq z'\}$. 
Returning to $k\geq3$, such a $d_k^-$-set is \emph{completed} if there exists a \emph{completing element} $z_k$ such that $\{w_k,\ldots,w_3;x,y;z_3,\ldots, z_{k-1},z_k\}$ is a $d_k$-interval. 
Two $d_k$-intervals are \emph{distinct} if their set symmetric difference is non-empty.
Let $k\geq 4$. Suppose $[w,z']$ is a $d_k^-$-interval  in which $u$ is the unique element covering $w$. 
If there exists $w'\neq w$ also covered by $u$ such that $[w',z']$ is also a  $d_k^-$-interval, then the $d_k^-$-intervals $[w,z']$ and $[w',z']$ \emph{overlap}. Overlapping $d_k^-$-intervals are shown in Figure~\ref{dk(1) and overlapping}b: They differ only in their minimal elements \cite{Okada}.

\vspace{.1in}
The following statement is the analog of Fact~\ref{get interval} for $k\geq 4$:
\vspace{.1in}
\begin{fact}\label{get interval general}
Let $k\geq 4$. Let $S=\{w_k,\ldots,w_3;x,y;z_3,\ldots, z_{k-1}\}$ be a $d_k^-$-set and let $z_k\in P$. 

\noindent(a) If $S\cup\{z_k\}$ is a $d_k$-interval, then $S\cup\{z_k\}=[w_k,z_k]$ and $z_k$ covers $z_{k-1}$.

\noindent(b) If $z_k$ covers $z_{k-1}$ and no other elements, then $S\cup\{z_k\}$ is the $d_k$-interval $[w_k,z_k]$

\end{fact}

\newpage

\noindent \textbf{Axioms} \  \textit{Let $k\geq3$. A poset satisfies the}

\vspace{.1in}

\noindent [\textit{Class I: Completion Axiom}]

\begin{axiom*}\hspace{-.075in}($Dk^-C$)  \hspace{.001in}
``$d_k^-$-sets are Completed" axiom if for each $d_k^-$-set $S$ there exists an element $z$ such that $S\cup \{z\}$ is a $d_k$-interval, the
\end{axiom*}

\vspace{.1in}

\noindent [\textit{Class II: Freeness Axiom}]

\begin{axiom*}\hspace{-.075in}(DkMF)  \hspace{.001in}
``$d_k$-interval Maxs are Free" axiom if the maximum element of each $d_k$-interval covers only (an) element(s) in that interval, the
\end{axiom*}

\vspace{.1in}

 \noindent [\textit{Classes I/ II: Completion/Freeness Axiom}]

\begin{axiom*}\hspace{-.075in}($Dk^-CF$)  \hspace{.001in}
``$d_k^-$-sets are Completed Freely" axiom if for each $d_k^-$-set $S$ there exists  an element $z$ such that $S\cup \{z\}$ is a $d_k$-interval with maximum element $z$ such that $z$ covers only (an) element(s) from $S$, the
\end{axiom*}

\vspace{.1in}

 \noindent [\textit{Class III: Forbidden Structure Axioms}]

  \begin{axiom*}\hspace{-.075in}($NODk^-$)  \hspace{.001in}
 ``No Overlapping $d_k^-$-sets" axiom if there do not exist overlapping $d_k^-$-sets, and the
 \end{axiom*} 

\begin{axiom*}\hspace{-.075in}(DkMD)  \hspace{.001in}
 ``$d_k$-interval Maxs are Distinct" axiom if the maximum elements of distinct $d_k$-intervals are distinct. \end{axiom*}

\begin{myremark}\label{obvious axiom implications kgeq3}
Each of these axioms subsumes the corresponding $k=3$ axiom, with NCC having been renamed NOD3$^-$. Axiom Dk$^-$CF obviously implies Dk$^-$C. Axioms Dk$^-$C and DkMF together imply Dk$^-$CF. 
\end{myremark}

\noindent \textbf{Properties}. \textit{Let $k\geq 3$. A poset has the}

\begin{property*}\hspace{-.075in}(YECOI)  \hspace{.001in}
 ``$Y$-stem Elements Covered Only Internally" property if each stem element of a $Y_k$-set is covered only by element(s) from that set, the
 \end{property*}
 
 \begin{property*}\hspace{-.075in}(N$\Lambda$Yk)  \hspace{.001in}
 ``No $\Lambda Y_k$-sets" property if there do not exists $\Lambda Y_k$-sets, and the \end{property*}

\begin{property*}\hspace{-.075in}(UCk)  \hspace{.001in}
``Unique Completion" property if each $d_k^-$-set has exactly one completing element.
\end{property*}


\section{Results for $\bold{\emph{k}\geq3}$}
\label{Propositions kgeq3}

Let $k\geq 3$ throughout this section.

\begin{myproposition}[\normalfont{YECOI}]\label{YECOI}  \hspace{.001in}
If a VT poset has No Triply Covereds, then it has $Y$-stem Elements Covered Only Internally.
\end{myproposition}

 \noindent Here is how we can extend a $d_k$-interval to a  $d_{k+1}^-$-interval:

\begin{mylemma}[\normalfont{Add To Tail $\#1$ (ATTk $\#1$)}]\label{ATTk1}  \hspace{.001in}
Consider any poset. Let $\{w_k,\ldots,w_3;x,y;z_3,\ldots,z_k\}$ be a $d_k$-interval and let $w_{k+1}$ be such that $w_{k+1}\rightarrow w_k$. If $[w_{k+1};x,y]$ is a $Y_{k+1}$-set and the neck elements of $[w_k,z_k]$ are free, then $[w_{k+1},z_k]$ is a $d_{k+1}^-$-interval.
\end{mylemma}

\noindent The hybrid axioms Dh$^-$CF for $3\leq h\leq k$ enable us to extend $Y_k$-sets to nice $d_k$-intervals:

\begin{mylemma}[\normalfont{$Y_k$-sets are Completed Freely (YkCF)}]\label{YkCF}  \hspace{.001in}
Suppose a poset is $Dh^-CF$ for $3\leq h\leq k$. If \\ $\{w_k,w_{k-1},\ldots,w_3;x,y\}$  is a $Y_k$-set, then there exist elements $z_3\rightarrow z_4\rightarrow\cdots\rightarrow z_k$ such that $[w_h,z_h]$ is a $d_h$-interval and $z_h$ is free for $3\leq h \leq k$.
\end{mylemma}

\noindent If in addition we assume the Class~III axiom at the next index, we can rule out a forbidden structure:

\begin{myproposition}[\normalfont N$\mathrm{\Lambda}$Yk]  \label{NlambdaYk} \hspace{.001in}
Suppose a poset is $NOD(k+1)^-$ and $Dh^-CF$ for $3\leq h\leq k$. Then there cannot exist $\Lambda Y_k$-sets.
\end{myproposition}

In the following theorem, Parts (\firstNODk), (\generalclassiccombo) (\secondNODk), (\firstDkMD), (\secondDkMD), and the second part of  (\UCkDkMF) respectively generalize Part (\NCCUC), Part (\classiccombo), the first part of Part (\secondNCCMaxsFree), the second part of Part (\MaxsDistinctMaxsFree), Part (\MaxsDistinct), and the first part of Part (\MaxsDistinctMaxsFree) of Theorem~\ref{k=3 locally finite theorem}. The first part of Part (\UCkDkMF)  generalizes the third part of Theorem~\ref{k=3 locally finite theorem}(\DAIUC), once D3$^-$C has been strengthened to Dh$^-$CF for $3\leq h \leq k$.

\begin{theorem} \label{kgeq3 locally finite theorem} The implications in Table~\ref{kgeq3 locally finite table} hold in a poset; here the letter ``h" indicates that the axiom is to be assumed for $3\leq h\leq k$. \\
\begin{center}
\begin{tabular}{ c @{\hspace{.25in}}c  c   @{\hspace{0in}}c   @{\hspace{.15in}}c  c c  c c  c c}

	 && \hspace{-.35in} \underline{I}  & & \hspace{-.05in}\underline{II} &  &\underline{III}  && \underline{Property} & &\underline{Conclusion(s)} \\
(\firstNODk) && \hspace{-.2in}VT, $Dk^-C$ &&---  && DkMD && NTC && $NODk^-$ \\
(\generalclassiccombo)  && & ${Dk^-CF}$ & &&  --- && UCk && $Dk^-C$ + DkMF \\
(\secondNODk)   &  &  & $Dk^-CF$ & & & DkMD && --- & & $NODk^-$ \\
(\firstDkMD)    &&  &$Dh^-CF$& && $NODh^-$  && --- && DkMD \\
(\secondDkMD)   &&  \hspace{-.3in}--- && DhMF&& $NODh^- $  && --- && DkMD \\
(\UCkDkMF)   &&  &$Dh^-CF$ &  && $NODk^-$ && --- && UCk, DkMF\\
\end{tabular}
\captionof{table}{Implications for Theorem~\ref{kgeq3 locally finite theorem}}
\label{kgeq3 locally finite table}
\end{center}
\end{theorem}

\noindent The hypotheses for part (\secondNODk) and (\firstDkMD) will be used in Section~\ref{d-complete definitions} to define $d_k$-complete and $d_h$-complete posets respectively.

The remaining results assume the No Triply Covereds property. In practice these results may most often be applied to finite posets via the following: By the second part of Theorem~\ref{k=3 finite theorem}(\thirdNCCsecondNTCUC)  (or Theorem~\ref{k=3 finite theorem}(\firstNTC)), a D3$^-$C poset has the No Triply Covereds property if it is finite and is D3MD (or is D3$^-$CF).

\begin{mylemma}[\normalfont{Add To Tail $\#2$ (ATTk $\#2$)}]\label{ATTk2}  \hspace{.001in}
Consider a VT poset that has No Triply Covereds. Let $\{w_k,\ldots,w_3;\\x,y;z_3,\ldots,z_k\}$ be a $d_k$-interval and let $w_{k+1}$ be such that $w_{k+1}\rightarrow w_k$. If $[w_{k+1};x,y]$ is a $Y_{k+1}$-set, then $[w_{k+1},z_k]$ is a $d_{k+1}^-$-interval.
\end{mylemma}

\begin{mylemma}[\normalfont{$Y_k$-sets are Completed  (YkC)}]\label{YkC}  \hspace{.001in}
Consider a poset that is $Dh^-C$ for $3\leq h\leq k$ and has No Triply Covereds. If  $\{w_k,w_{k-1},\ldots,w_3;x,y\}$  is a $Y_k$-set, then there exist elements $z_3\rightarrow z_4\rightarrow\cdots\rightarrow z_k$ such that $[w_h,z_h]$ is a $d_h$-interval for $3\leq h \leq k$.
\end{mylemma}

\noindent The next result, which obtains the unique completion property a second time, generalizes the result of following Theorem~\ref{k=3 locally finite theorem}(\NCCUC) by the last part of Theorem~\ref{k=3 locally finite theorem}(\DAIUC):

\begin{myproposition}[\normalfont{UCk $\#2$}]\label{UCk}  \hspace{.001in} Consider a  poset that is  $Dh^-C$ for $3\leq h\leq k$ and has No Triply Covereds. If it is DkMD, then it has Unique Completion at $k$.
\end{myproposition}


\section{Proofs for $\bold{\emph{k}\geq3}$}
\label{proofs kgeq3}

\begin{proof}\hspace{-.035in}\textit{of Proposition~\ref{YECOI}}.
Suppose some stem element $w_i$ of a $Y_k$-set $[w_k;x,y]$ is covered by some $u\not\in[w_k;x,y]$. Apply UPUE from Theorem~\ref{k=3 locally finite theorem}(\UPUE) to $w_i\leq w_3$, $w_i\rightarrow u$, and $u\not\leq w_3$ to produce $z\not\in\{x,y\}$ such that $w_3\rightarrow z$.
\qed
\end{proof} 

\begin{proof}\hspace{-.035in}\textit{of Lemma~\ref{ATTk1}}.
Let $u\in[w_{k+1},z_k]$ be maximal such that $u\not\in\{w_{k+1}\}\cup[w_k,z_k]$. Here maximality implies that $u$ is covered by some element $v$ of $[w_k,z_k]$. This $v$ cannot be a neck element, since those are free. But $v\in\{w_k,\ldots,w_3,x,y\}$ with $u\geq w_{k+1}$ would violate $[w_{k+1};x,y]$ being a $Y_{k+1}$-set. So there is no such $u$. Hence $[w_{k+1},z_k]\cong {\text{\d}_{k+1}(1)}-\{f_{k+1}\}$.
\qed
\end{proof}

\begin{proof}\hspace{-.035in}\textit{of Lemma~\ref{YkCF}}.
Here D3$^-$CF says that $\{w_3;x,y\}$ is freely completed with a $z_3$. And $\{w_4,w_3;x,y\}$ is a $Y_4$-set. So Lemma ATT3$\#1$  says that $[w_4,z_3]$ is a $d_4^-$-interval. Now repeatedly alternate the application of Dh$^-$CF and then ATTh$\#1$ for $4\leq h< k$ to construct the $d_h^-$-intervals $[w_{h+1},z_{h}]$. Finish with Dk$^-$CF.
\qed
\end{proof}

\begin{proof}\hspace{-.035in}\textit{of Proposition~\ref{NlambdaYk}}.
Suppose $[u,v;x,y]$ is a $\Lambda Y_k$-set with $\{u,v\}\rightarrow w_k\in[u,v;x,y]$. Here $[w_k;x,y]$ is a $Y_k$-set. Apply Lemma YkCF to construct a $d_k$-interval $[w_k,z_k]$ with free neck elements $z_3,\ldots,z_k$. Note that $[u;x,y]$ and $[v;x,y]$ are $Y_{k+1}$-sets. So Lemma ATTk$\#1$ says that $[u,z_k]$ and $[v,z_k]$ are $d_{k+1}^-$-intervals. But they are overlapping.
\qed
\end{proof}

\begin{proof}\hspace{-.035in}\textit{of Theorem~\ref{kgeq3 locally finite theorem}}.
Part (b) follows from the definitions for $k\geq 3$. Since it has been noted how the other parts reduce to parts of Theorem~\ref{k=3 locally finite theorem} when $k=3$, suppose $k\geq4$.

(\firstNODk) Let $\{w_k, w_{k-1},\ldots,w_3; x,y;z_3,\ldots,z_{k-1}\}$ and $\{w_k',w_{k-1},\ldots,w_3;x,y;z_3,\ldots,z_{k-1}\}$ be overlapping $d_k^-$-intervals. Apply Dk$^-$C to obtain $z_k$ such that $[w_k,z_k]$ is a $d_k$-interval. By DkMD, there exists $u\in[w_k',z_k]$ with $u\not\in\{w_k,w_k'\}\cup[w_{k-1},z_k]$. Let $v$ be the minimal such $u$. Since $w_{k-1}<v$ is not possible, it must be that $w_k'\rightarrow v$.  Since $[w_k';x,y]$ is a $Y_k$-set, this violates Proposition YECOI.

(\secondNODk)
Let $[w_k,z_{k-1}]$ and $[w_k',z_{k-1}]$ be overlapping $d_k^-$-intervals. Apply Dk$^-$CF to obtain $z_k$ such that $[w_k,z_k]$ is a $d_k$-interval with $z_k$ free. Since $z_k$ covers only $z_{k-1}$, we see that $[w_k',z_k]$ is a $d_k$-interval. This contradicts DkMD.

(\firstDkMD)
Let $\{w_k,\ldots,w_3;x,y;z_3,\ldots,z_k\}$ and $\{a_k,\ldots,a_3;b,c;f_3,\ldots, f_k\}$ be two $d_k$-intervals with $z_k=f_k$. Apply Dk$^-$CF to $[w_k,z_{k-1}]$ to produce $z_k'$ such that $[w_k,z_k']$ is a $d_k$-interval with $z_k'$ free. Suppose $z_k'\neq z_k$. Here $[x,y;z_k,z_k']$ is a $\Lambda Y_{k-1}$-set. This would contradict Proposition N$\mathrm{\Lambda}$Y(k-1), and so $z_k'=z_k$. Hence $z_k$ is free, which implies $z_{k-1}=f_{k-1}$. This argument can be repeated to conclude that $z_h=f_h$ for $k\geq h\geq 3$. From Theorem~\ref{k=3 locally finite theorem}(\MaxsDistinctMaxsFree) we have D3MD. So $\{x,y\}=\{b,c\}$ and $w_3=a_3$. By NODh$^-$ for $4\leq h\leq k$ we have $w_h=a_h$. Hence $[w_k,z_k]=[a_k,f_h]$.

(\secondDkMD)
Let $\{w_k,\ldots,w_3;x,y;z_3,\ldots,z_k\}$ and $\{a_k,\ldots,a_3;b,c;f_3,\ldots,f_k\}$ be two $d_k$-intervals with $z_k=f_k$. By DhMF for $k\geq h\geq 4$ we have $z_{h-1}=f_{h-1}$. By D3MF we have $\{x,y\}=\{b,c\}$. Now NCC requires $w_3=a_3$. Finish as in the proof of Part (\firstDkMD).

(\UCkDkMF)
We note that for UCk we will not need the freeness of Dh$^-$CF at $h=k$, but only for $3\leq h\leq k-1$.
Let $\{w_k,\ldots ,w_3;x,y;z_3,\ldots, z_{k-1}\}$ be a $d_k^-$-interval. Then  Dk$^-$C gives a $z_k$ so that $[w_k,z_k]$ is a $d_k$-interval. Suppose that $z_k'\neq z_k$ also completes $[w_k,z_{k-1}]$, now to a $d_k$-interval $[w_k,z_k']$. Here $[x,y;z_k,z_k']$ is a $\Lambda Y_{k-1}$-set. Since this would contradict Proposition N$\mathrm{\Lambda}$Y(k-1), Property UCk holds. For DkMF, let $\{w_k,\ldots,w_3;x,y;z_3,\ldots, z_k\}$ be a $d_k$-interval. Then Dk$^-$CF gives a free completing element $z_k'$ for the $d_k^-$-interval $[w_k,z_{k-1}]$. Here UCk implies that $z_k'=z_k$, and so $z_k$ is free.
\qed
\end{proof}

\begin{proof}\hspace{-.035in}\textit{of Lemma~\ref{ATTk2}}.
Let $u\in[w_{k+1},z_k]$ be minimal such that $u\not\in\{w_{k+1}\}\cup[w_k,z_k]$. Here minimality implies that $u$ covers some element $v$ of $\{w_{k+1}\}\cup[w_k,z_{k-1}]$. Since $u\leq z_k$, we cannot have $v\in[w_k,z_{k-1}]$. So $v=w_{k+1}$, and $w_{k+1}\rightarrow u$ violates Proposition YECOI.
\qed
\end{proof}

\begin{proof}\hspace{-.035in}\textit{of Lemma~\ref{YkC}}.
In the proof of Lemma~\ref{YkCF}, replace each instance of `Dh$^-$CF' with `Dh$^-$C', replace each instance of `ATTh$\#1$' with `ATTh$\#2$', and delete all references to `free' completions.
\qed
\end{proof}

\begin{proof}\hspace{-.035in}\textit{of Proposition~\ref{UCk}}.
Given the generalization remark, suppose $k\geq4$. Let $\{w_k,\ldots,w_3;x,y;z_3,\ldots,z_{k-1}\}$ be a $d_k^-$-interval. Then Dk$^-$C gives a $z_k$ so that $[w_k,z_k]$ is a $d_k$-interval. Suppose that $z_k'\neq z_k$ also completes $[w_k,z_{k-1}]$, now to a $d_k$-interval $[w_k,z_k']$. Here $[z_3;z_k,z_k']$ is a $Y_{k-1}$-set.  Lemma Y(k-1)C constructs $u_{k-1}\rightarrow u_{k-2}\rightarrow\cdots\rightarrow u_3$ such that $[z_h,u_h]$ is a $d_{k+2-h}$-interval for $k-1\geq h \geq 3$. So $[z_3,u_3]$ is a $d_{k-1}$-interval. Here $[x;z_k,z_k']$ and $[y;z_k,z_k']$ are $Y_k$-sets. Use Lemma ATT(k-1)$\#2$ on each of them to see that $[x,u_3]$ and $[y,u_3]$ are $d_k^-$-intervals. They are overlapping, which contradicts Theorem~\ref{kgeq3 locally finite theorem}(\firstNODk). 
\qed
\end{proof}


\section{Definitions for $\bold{\emph{d}}$-complete posets and their equivalences}\label{d-complete definitions}

We continue to consider locally finite posets.
From \cite{Japan}, here are our currently preferred (and shortest)
definitions of $d_k$-complete, $d_{\leq k}$-complete, and $d$-complete posets; the parenthetical words are to be invoked when $k=3$:

\begin{mydefinition}[\normalfont{K\^{o}ky\^{u}roku}]\label{kk}
A poset is $d_k$-complete if for every $d_k^-$-set $S$ there exists an element that covers exactly the maximal element(s) of $S$ and that does not cover (both of) the maximal element(s) of any other $d_k^-$-set.
It is $d_{\leq k}$-complete if it is $d_h$-complete for every
$3\leq h\leq k$ and it is $d$-complete if it is $d_k$-complete for every $k\geq 3$.
\end{mydefinition}

 Alternatively, one could define these three notions using any one of the four combinations of axioms from Sections~\ref{k=3 propositions} and~\ref{Propositions kgeq3} that are presented in Table~\ref{d-complete definitions table}.
Once the remark in Section~\ref{k=3 propositions} concerning the inadvertent double-defining of
$d_3$-complete in \cite{DDCT} \cite{Wave} is taken into account,  Combination (a) of Table~\ref{d-complete definitions table} at $h=k$ was essentially used in those papers to define $d_k$-complete and $d$-complete finite posets. A nicely worded version of that combination appeared in \cite{Okada}. We refer to it as the \textit{Classic} definition.

Now we relate the (K\^{o}ky\^{u}roku) definitions of $d_k$-complete and $d_{\leq k}$-complete to the four combinations of axioms. For the notion of $d_k$-complete, some of the $5\times 4=20$ possible implications among these $4+1=5$ criteria do not hold true at $h=k$ alone. But some do hold true. See Section~\ref{equivalence and general properties proofs} for further remarks. For the notion of $d_{\leq k}$-complete, all of these 20 implications hold true when one or more of the two or three hypothesis axioms is assumed for $3\leq h\leq k$:

\begin{theorem}\label{equivalence theorem}
 Let $k\geq3$.  A poset is $d_{\leq k}$-complete
if and only if it satisfies any one of the Combinations (a) - (d) of axioms displayed in Table~\ref{d-complete definitions table} for $3\leq h\leq k$.  Hence it is $d$-complete if and only if it satisfies any one of these
combinations of axioms for $k\geq 3$.

\begin{center}
\begin{tabular}{@{\hspace{.2in}}c   @{\hspace{.4in}}c     @{\hspace{-.1in}}c  c   @{\hspace{-.1in}}c   @{\hspace{.1in}}c  @{\hspace{.3in}}c c @{\hspace{.8in}} c}
	& &  && \hspace{-.1in} \underline{I}  & & \hspace{-.05in}\underline{II} &  &\underline{III}  \\
(a)  &&  && $Dh^-C$     &&       DhMF &&   $NODh^-$ \\
(b)    && &&  $Dh^-C$     &&       DhMF   && DhMD\\
(c)        &&&    && \hspace{.2in}$Dh^-CF$       &   &&$NODh^-$  \\
(d)     &&       && &  \hspace{.2in}$Dh^-CF$     &&&      DhMD
\end{tabular}
\captionof{table}{Combinations (a) - (d) of axioms for $h\geq 3$}\label{d-complete definitions table}
\end{center}

\end{theorem}

\section{Properties of $\bold{\emph{d}}$-complete posets}\label{properties of d-complete posets}

We take note of two important facts that are not used in this paper:

\begin{fact}
Any ``filter" of a $d$-complete poset is $d$-complete, as is the disjoint union of two $d$-complete posets. 
\end{fact}

\noindent For the first statement, note that removing an ``ideal" of elements to produce a filter of $P$ does not adversely affect the satisfaction of the $d$-complete requirements for $P$.

The following theorem applies the results of Sections~\ref{k=3 propositions} and~\ref{Propositions kgeq3} to $d$-complete posets:

\begin{theorem}\label{all axioms and props} Let $P$ be a $d$-complete poset.

 \noindent (a) The poset $P$ satisfies all of the axioms defined in Sections~\ref{k=3 axioms and properties} and~\ref{kgeq3 axioms and properties} and possesses the following properties defined there: UPUE, DAI, UT, and UCk.

\noindent (b) If $P$ is finite, it also possesses all of the properties defined in Sections~\ref{k=3 axioms and properties} and~\ref{kgeq3 axioms and properties}, except for UM and CLMEE when $P$ is not connected.
\end{theorem}

\noindent Here Part (a) follows from the observation that by Theorem~\ref{equivalence theorem} the hypotheses of all of the Propositions and Theorems in Sections~\ref{k=3 propositions} and~\ref{Propositions kgeq3}  that do not assume finiteness, NTC, or SS are satisfied.
 For Part (b), note that if $P$ is finite, then the hypotheses of all of the Propositions and Theorems in Sections~\ref{k=3 propositions} and~\ref{Propositions kgeq3} are satisfied, apart from Theorem~\ref{k=3 finite theorem}(\UMCLMEE). 

We now study the necks and tails of DTD intervals in $d$-complete posets. In the following statement, Part (a) restates Axiom DhMF for $3\leq h \leq k$ and Part (b) follows immediately from  Proposition~YECOI of Section~\ref{Propositions kgeq3}:

\begin{fact}\label{neck/tail fact}
Let $P$ be a $d$-complete poset. Let $k\geq 3$. Let $\{w_k,\ldots, w_3;x,y;z_3,\ldots,z_k\}$ be a $d_k$-interval. 

\noindent(a) A neck element $z_i$ of $[w_k,z_k]$ can cover only element(s) in $[w_k,z_k]$: If $z_i$ is strict (i.e. $i\geq 4$), then it covers only $z_{i-1}$. Otherwise (i.e. $i=3$) it covers only $x$ and $y$.

\noindent(b) If $P$ has the No Triply Covereds property, then a tail element $w_i$ of $[w_k,z_k]$ can be covered by only element(s) in $[w_k,z_k]$: If $w_i$ is strict (i.e. $i\geq4$), then it is covered only by $w_{i-1}$. Otherwise (i.e. $i=3$) it is covered only by $x$ and $y$.
\end{fact}

\noindent The next result says that the necks (tails) of two DTD intervals may intersect only in a particular way. For Part (b), keep in mind that every finite $d$-complete poset has the NTC property. 

\begin{myproposition} \label{neck/tail intersection}
Let $P$ be a $d$-complete poset. Let $3\leq k\leq k'$. 

\noindent(a) If there exists an element that is both a neck element for a $d_k$-interval $[w_k,z_k]$ and a neck element for a $d_{k'}$-interval $[a_{k'},f_{k'}]$, then $[w_k,z_k]\subseteq[a_{k'},f_{k'}]$.

\noindent(b) If $P$ has the No Triply Covereds property and there exists an element that is both a tail element for a $d_k$-interval $[w_k,z_k]$ and a tail element for a $d_{k'}$-interval $[a_{k'},f_{k'}]$, then $[w_k,z_k]\subseteq[a_{k'},f_{k'}]$.
\end{myproposition}

\begin{mycorollary}\label{overlap corollary}
Let $P$ be a $d$-complete poset. Let $k\geq 3$. Let $\{w_k,\ldots,w_3;x,y;z_3,\ldots,z_k\}$ be a $d_k$-interval. 

\noindent For each $k'\geq k$ there is at most one $d_{k'}$-interval that contains $[w_k,z_k]$, and if such an interval exists it must be of the form $\{w_{k'},\ldots,w_k,\ldots w_3;x,y;z_3,\ldots,z_k,\ldots, z_{k'}\}$. 
 A neck element $z_i$ of $[w_k,z_k]$ can be a neck element of only those $d_{k'}$-intervals  $\{w_{k'},\ldots,w_k,\ldots,w_3;x,y;z_3,\ldots,z_k,\ldots,z_{k'}\}$.
 If $P$ has the No Triply Covereds property, then a tail element $w_i$ of $[w_k,z_k]$ can be a tail element of only those $d_{k'}$-intervals $\{w_{k'},\ldots,w_k,\ldots,w_3;x,y;\\z_3,\ldots,z_k,\ldots,z_{k'}\}$.

\end{mycorollary}


\section{Proofs of equivalences and properties }\label{equivalence and general properties proofs}

\begin{proof} \hspace{-.035in}   \textit{of Theorem~\ref{equivalence theorem}}.
Table~\ref{proof of equiv table} presents six implications for parts of the $d_{\leq k}$-complete statement. The first five come from Sections~\ref{kgeq3 axioms and properties} and~\ref{Propositions kgeq3}. Here the `$h=k$' and `$h\leq k$' entries under ``Realm" indicate whether the hypothesis of the implication needs to assume that the axioms at hand hold merely at $k$ or it needs to assume that the axioms hold for all $3\leq h\leq k$. The last implication is verified by composing two of the earlier implications and then remembering one of its hypotheses. So we can finish this proof by relating the K\^{o}ky\^{u}roku definition to any of these combinations of axioms. We show that Combination (d) $\Rightarrow$ K\^{o}ky\^{u}roku and that K\^{o}ky\^{u}roku $\Rightarrow$ Combination (c), both within the $h=k$ realm. Let $k\geq 3$.

\begin{center}
\begin{tabular}{c @{\hspace{.6in}}c @{\hspace{.6in}}c}
\underline{Implication} & \underline{Realm} & \underline{Citation/Justification} \\
(a) $\Rightarrow$ (c) & $h=k$ & Remark~\ref{obvious axiom implications kgeq3} \\
(c) $\Rightarrow$ (a) & $h\leq k$ & Remark~\ref{obvious axiom implications kgeq3} + Theorem~\ref{kgeq3 locally finite theorem}(\UCkDkMF) \\
(d) $\Rightarrow$ (c) & $h=k$ & Theorem~\ref{kgeq3 locally finite theorem}(\secondNODk) \\
(c) $\Rightarrow$ (d) & $h\leq k$ & Theorem~\ref{kgeq3 locally finite theorem}(\firstDkMD) \\
(b) $\Rightarrow$ (d) & $h=k$ & Remark~\ref{obvious axiom implications kgeq3} \\
(a) $\Rightarrow$ (b) & $h\leq k$ & (a) $\Rightarrow$ (c) $\Rightarrow$ (d); + (a)
\end{tabular}
\captionof{table}{Implications in Theorem~\ref{equivalence theorem}}\label{proof of equiv table}
\end{center}

Suppose Combination (d) holds at $h=k$. Let $S$ be a $d_k^-$-set with minimum element $w_k$ and maximum element $z_{k-1}$ (when $k\geq 4$) or maximum elements $\{x,y\}$ (when $k=3$). Then Dk$^-$CF gives some $z_k$ such that $S\cup \{z_k\}$ is the $d_k$-interval $[w_k,z_k]$ and $z_k$ covers only element(s) from $S$. Facts~\ref{get interval general}(a) and \ref{get interval}(a) imply that $z_k$ covers exactly these maximum element(s) of $S$.
Suppose that $z_k$ covers the maximum element(s) $f_{k-1}$ (or $\{b,c\}$) of some $d_k^-$-set $T$, whose minimum element is denoted $a_k$. Since $z_k$ covers $z_{k-1}$ (or $\{x,y\}$) exactly, we have $z_{k-1}=f_{k-1}$ or $\{x,y\}=\{b,c\}$. Since $z_k$ covers $f_{k-1}$ (or $\{b,c\}$) exactly, Facts~\ref{get interval general}(b) and \ref{get interval}(b) say that $T\cup\{z_k\}$ is the $d_k$-interval $[a_k, z_k]$.  Here DkMD says that $[w_k, z_k]=[a_k, z_k]$, and so $T$ must coincide with $S$. Hence K\^{o}ky\^{u}roku holds.

Suppose the K\^{o}ky\^{u}roku definition of $d_k$-complete holds at $k\geq 3$. Let $S$ be a $d_k^-$-set with minimum element $w_k$ and maximum element $z_{k-1}$ (when $k\geq 4)$ or maximum elements $\{x,y\}$ (when $k=3$). Then by K\^{o}ky\^{u}roku there exists some $z_k$ that covers these maximum element(s) of $S$ exactly and that does not cover the maximum element(s) of any other $d_k^-$-set. Facts~\ref{get interval general}(b) and \ref{get interval}(b) say that $S\cup\{z_k\}$ is the $d_k$-interval $[w_k, z_k]$, and so Dk$^-$CF is satisfied.
Suppose $\{w_k, w_{k-1}, \ldots, w_3;x,y;z_3,\ldots, z_{k-1}\}$ and $\{w_{k'},w_{k-1},\ldots, w_3;x,y;\\z_3,\ldots, z_{k-1}\}$ are overlapping $d_k^-$-sets. There exists some $z_k$ that covers $z_{k-1}$ (or $\{x,y\}$) exactly and that does not cover the maximal element(s) of any other $d_k^-$-set. Since $w_{k'}\neq w_k$ and the maximal element(s) of $\{w_{k'}, w_{k-1}, \ldots, w_3; x,y; z_3,\ldots, z_{k-1}\}$ is (are) $z_{k-1}$ (or $\{x,y\}$), we see that $z_k$ covers the maximal element(s) of another $d_k^-$-set. This contradiction implies that NODk$^-$ holds. So Combination (c) is satisfied at $h=k$.
\qed
\end{proof}

\begin{proof}\hspace{-.035in}\textit{of Proposition~\ref{neck/tail intersection}}.
Let $3\leq k\leq k'.$ Let $\{w_k,\ldots,w_3;x,y;z_3,\ldots,z_k\}$ be a $d_k$-interval and $\{a_{k'},\ldots,a_3; \newline b,c;f_3,\ldots,f_{k'}\}$ be a $d_{k'}$-interval. 

(a) Let $3\leq j\leq k'$ be maximal such that there exists $3\leq i\leq k$ with $z_i=f_j$. Let $m\defeq \text{min}\{i-3,j-3\}$.
Since $[w_k,z_k]$ is a $d_k$-interval, if $m\geq 1$ use Fact~\ref{neck/tail fact}(a) twice to get  $z_{i-1}=f_{j-1}$. Similarly, if $m\geq 2$ then $z_{i-2}=f_{j-2}$.
Continuing downward, we conclude that $z_{i-l}=f_{j-l}$ for $l\in\{0,\ldots,m\}$. 
This includes the case $m=0$.
If $i\neq j$, then $\{i-m,j-m\}=\{3,h\}$ with $h> 3$. So one of $z_{i-m}$ and $f_{j-m}$ is a diamond top and thus covers two distinct elements while the other is the maximum element of a $d_h$-interval and thus by Fact~\ref{neck/tail fact}(a) can only cover one element. This contradicts $z_{i-m}=f_{j-m}$. Thus $i=j$ and $z_{i-l}=f_{i-l}$ for $l\in\{0,\ldots, i-3\}$.
Since $z_3=f_3$, Axiom D3MD implies that $[w_3,z_3]=[a_3,f_3]$. So $\{x,y\}=\{b,c\}$ and $w_3=a_3$.

Suppose $i<k$. Then the element  $z_{i+1}$ exists and the choice of $j$ implies that $z_{i+1}\neq f_{i+1}$. 
Since $[w_k, z_k]$ and $[a_{k'}, f_{k'}]$ are DTD intervals, we see that $[x,y;z_{i+1},f_{i+1}]$ is a $\Lambda Y_i$-set. This contradicts Proposition N$\mathrm{\Lambda}$Yi. Thus it must be that $i=k$. So $z_t=f_t$ for $t\in\{3,\ldots,k\}$.
Here $[w_k,z_k]$ and $[a_k,f_k]$ are both $d_k$-intervals with $z_k=f_k$. Hence DkMD implies $[w_k,z_k]=[a_k,f_k]$.
Therefore $[w_k,z_k]\subseteq[a_{k'},f_{k'}]$.

 (b) Let $3\leq j \leq k'$ be such that there exists $3\leq i\leq k$ with $w_i=a_j$. Let $m\defeq \text{min}\{i-3,j-3\}$. 
 The argument above can be ``reflected" to move up toward a diamond: Use Fact~\ref{neck/tail fact}(b) instead of Fact~\ref{neck/tail fact}(a) to soon obtain $w_{i-1}=a_{j-1}$. After analogizing five more sentences (again using Fact~\ref{neck/tail fact}(b) instead of Fact~\ref{neck/tail fact}(a)), we arrive at contradicting $w_{i-m}=a_{j-m}$. 
Thus $i=j$ and $w_{i-l}=a_{i-l}$ for $l\in\{0,\ldots,i-3\}$.
Since $w_3=a_3$, to avoid contradicting NTC it must be that $\{x,y\}=\{b,c\}$. Then $z_3=f_3$ to avoid contradicting NCC. Now note that $z_3$ is both a neck element for the $d_k$-interval $[w_k,z_k]$ and a neck element for the $d_{k'}$-interval $[a_{k'},f_{k'}]$. Then Part (a) implies that $[w_k,z_k]\subseteq [a_{k'},f_{k'}]$. \qed
\end{proof}

\begin{proof} \hspace{-.035in}\textit{of Corollary~\ref{overlap corollary}}.
Let $k'\geq k\geq3$ and let $\{w_k,\ldots,w_3;x,y;z_3,\ldots,z_k\}$ be a $d_k$-interval. 
Suppose two $d_{k'}$-intervals contain $[w_k,z_k]$. Their elbows must pairwise coincide with $\{x,y\}$ and $z_k$ is a neck element of both $d_{k'}$-intervals. Using Proposition~\ref{neck/tail intersection}(a) for two containments, we find that the two $d_{k'}$-intervals must be equal and have the claimed form. 
Let $z_i$ (or $w_i$) be a neck (respectively tail) element of $[w_k,z_k]$. If $z_i$ (or $w_i$) is also a neck (respectively tail) element of a $d_{k'}$-interval, then using Proposition~\ref{neck/tail intersection} and the first statement we find that that $d_{k'}$-interval must be the unique $d_{k'}$-interval that contains $[w_k,z_k]$. \qed
\end{proof}


\section{Other work on $\bold{\emph{d}}$-complete posets}\label{other work on d-complete posets}\label{other work}

For the most part,  we list only papers that work in a substantive fashion with  $d$-complete posets that are more general than filters of minuscule posets or rooted trees.   We include structures that are closely related to $d$-complete posets:  $\lambda$-minuscule elements of Kac-Moody Weyl groups, their heaps, and Nakada's ``generalized Young diagrams".

	Colors play no role in some appearances of  $d$-complete posets,  beginning with their classification \cite{DDCT} and continuing with the jeu de taquin result of \cite{JDT}.   Ishikawa and Tagawa used determinants and Pfaffians \cite{IT1} to prove that many classes of slant irreducible  $d$-complete posets possess Stanley's hook product property.   For standard shifted Young tableaux Konvalinka gave \cite{Kon} a bijective proof of the branching recursion that implies the hook product enumeration formula,  and he began to develop this approach for proving the hook product formula for counting linear extensions of  $d$-complete posets.   Riegler and Neumann \cite{RiNe} use jeu de taquin slides with respect to a fixed linear extension of a poset  $P$  to sort any labelling into an linear extension of  $P$.   They began to study this for  $d$-complete posets,  showing that the linear extensions produced are uniformly distributed when  $P$  is a filter of  \d$_n(1)$  (and hence  $d$-complete),  but not when  $P$  is a non-chain proper ideal  of  \d$_n(1)$  (and hence not  $d$-complete).   The website \cite{GaPr} has lists of connected  $d$-complete posets with up to 9 elements and a Mathematica procedure that determines whether a poset is  $d$-complete.

	Some appearances of  $d$-complete posets have initial statements that refer to uncolored structures from pre-existing combinatorial problems,  but at the same time have fuller colored statements or have proofs that refer to a colored version of the  $d$-complete poset.   In addition to the hook product identity \cite{Japan} found with Peterson for  $d$-complete posets,  this remark also applies to the generalizations of that identity found by Ishikawa and Tagawa for leaf posets \cite{IsTa}.   Okamura referred to the classification of  $d$-complete posets to give a case-by-case probabalistic proof \cite{Okam} of the hook product formula for counting the number of linear extensions of a  $d$-complete poset.   Nakada's results concern ``generalized Young diagram" posets that are formed from Kac-Moody roots:   In \cite{Nak1} his fractional ``colored hook formula" was a multivariate generalization of the formula used by Greene, Nijenhuis, and Wilf for their probabalistic proof of the hook product formula for counting standard Young tableaux.   In \cite{Nak2} he presented his version of the multivariate hook product identity of \cite{Japan}.   Nakada and Okamura noted \cite{NaOk} that a uniform probability algorithm proof of a product formula for counting linear extensions of these posets that was analogous to that of \cite{Okam} could be deduced in this context from \cite{Nak1}.   After proving  $(q,t)$-generalizations of multivariate hook product identities for reverse plane partitions on shapes and shifted shapes,  Okada conjectured \cite{Okada} an extension of it that would  $(q,t)$-generalize the hook product identity of \cite{Japan} for $d$-complete posets.   He confirmed this for rooted trees,  and Ishikawa confirmed \cite{Ishi} it for two of the simpler classes of slant irreducible  $d$-complete posets.   Kawanaka extended \cite{Kaw1} the Sato-Welter-Sprague-Grundy winning strategy for nim from shapes to  $d$-complete posets.   Later he introduced \cite{Kaw2} ``finitely branching principal plain algorithm" games;  it can be seen using \cite{Wave} that portions of the digraphs of these games arise in his Theorem~6.3 as the Hasse diagrams of lattices of filters of  $d$-complete posets.   Uncolored  $d$-complete posets can serve as ``boards" on which jeu de taquin rectification procedures are performed during the computation of cohomology products for some Schubert varieties in some flag varieties.   The ``squares" of these boards do not need to be colored for the sliding mechanics,  but they need to be colored when one labels the Schubert varieties with elements of the Weyl group.    Some such results of Chaput and Perrin \cite{ChPe} for Kac-Moody flag varieties use the well defined jeu de taquin rectification result of \cite{JDT} for some  $d$-complete posets that are not filters of minuscule posets.   The  $K$-theoretic Littlewood-Richardson results of Buch and Samuel \cite{BuSa} refer only to minuscule posets,  as do several cohomology computation references of \cite{BuSa}.

	Colors play a central role in some appearances of  $d$-complete posets,  beginning with their first formulation in \cite{Mar}.   Earlier,  the product formula on p.~348 of \cite{BLP} for the number of linear extensions of a minuscule poset did not refer to colors.   However,  Theorem~11 there described a minuscule poset as a poset of certain colored coroots for its associated Weyl group.   Combining the remark on pp.~345-346 with Theorem~11,  in hindsight that product formula also expressed the number of reduced decompositions of a minuscule element of a finite Weyl group (a colored problem) as a product over a poset of colored coroots.   Peterson extended this product-over-roots formula \cite{Carrell} for reduced decompositions to  $\lambda$-minuscule elements of Kac-Moody Weyl groups.   For further information on Peterson's work and the development of the notion of  $d$-complete in \cite{Mar} from the work in \cite{BLP},  which later led to \cite{Japan},  see Section~13 of \cite{Japan}.   Nakada's overview \cite{Nak3} of \cite{Nak1}, \cite{Nak2}, and \cite{NaOk} notes that Peterson's formula can be deduced from the main result of any of those papers.   Given the connection between the linear extensions of colored  $d$-complete posets and such reduced decompositions that was described in \cite{Wave} for the simply laced cases,  a closely related hook product formula for this number can be deduced from \cite{Okam} or \cite{Japan}.   Also via this connection,  the classification of  $d$-complete posets in \cite{DDCT} gave a classification of the  $\lambda$-minuscule elements of simply laced Kac-Moody Weyl groups.   Stembridge extended \cite{Ste} this classification to all symmetrizable Kac-Moody Weyl groups.   There Theorem~5.5 extended Theorem~11 of \cite{BLP} to use posets of coroots to describe the heaps of the  $\lambda$-minuscule elements in all symmetrizable Kac-Moody Weyl groups.   Kleshchev's and Ram's Theorem~3.10 of \cite{KlRa} can be seen to be saying that the dimensions of certain homogenous irreducible modules of Khovanov-Lauda-Rouquier algebras are equal to the number of linear extensions of associated  $d$-complete posets.   When the hook product expression of \cite{Okam} or \cite{Japan} is applied here,   this theorem generalizes the fact that the dimensions of the irreducible representations of the symmetric group are given by the FRT hook formula for enumerating standard Young tableaux.

	Green's ``full heaps" \cite{Gre} are candidates to be regarded as locally finite colored  $d$-complete posets once that definition is finalized;  they play a central role in that book.   Our Figure~\ref{e7(7) E6}b appears as the full heap of his Figure~6.13.   Lax refers to several of the axioms for  $d$-complete and colored  $d$-complete posets when he uses minuscule posets to give uniform derivations \cite{Lax} of the ``extreme" Pl\"{u}cker relations for the embeddings of minuscule flag varieties. Michael Strayer has shown (personal communication) that a finite poset $P$ can be colored in such a way that the lattice $J(P)$ carries a representation of a simply laced Kac-Moody Borel derived subalgebra in a certain natural ``minuscule" fashion exactly when $P$ is a simply colored $d$-complete poset.   

\begin{Added Notes}
Before this paper, the notion of ``$d$-complete" was considered only for finite posets.  As this paper was being written, it was observed that most of the axioms and definitions for finite $d$-complete posets continued to work well for locally finite posets,  without additions or modifications.  However,  recent work by Michael Strayer and the first author indicates that it will be useful in the future to require that the No Triply Covereds property holds as an axiom for an infinite locally finite poset to be called $d$-complete. Then the poset in Figure~\ref{infinite d-complete posets}a will no longer qualify to be $d$-complete.
            

Kim and Yoo evaluate integrals of the $q$-Selberg kind to give a new (class-by-class) proof \cite{KiYo} of Stanley's hook product property for all $d$-complete posets.  Naruse and Okada use formulas for products in the equivariant $K$-theory of Kac-Moody flag varieties to re-prove and generalize \cite{NrOk} the multivariate hook product identity of \cite{Japan} for $d$-complete posets.
\end{Added Notes}

\begin{acknowledgements}
We thank Soichi Okada, Alexander Kleshchev, and Arun Ram for their help on Section~\ref{other work}, and Michael Strayer for helpful remarks on the exposition.
\end{acknowledgements}

\end{document}